\documentclass[10pt,twoside]{article}
\usepackage[margin=1in]{geometry}
\usepackage{latexsym,amssymb,amsmath,amsfonts,xcolor,cases,mathtools}
\usepackage{graphicx,amsthm}
\usepackage{epstopdf,chngcntr}
\usepackage{natbib}
\usepackage{hyperref}
\usepackage{fancyhdr}

\numberwithin{equation}{section}

\hypersetup{
	colorlinks = true,
	linkcolor = {blue},
	citecolor = {blue},
	urlcolor  = {blue},}

\newcommand{\RE}[1]{\textup{Re}\left\{#1\right\}}
\newcommand{\IM}[1]{\textup{Im}\left\{#1\right\}}
\newcommand{\CC}{{\mathbb{C}}}

\newcommand{\ZZ}{{\mathbb{Z}}}

\newcommand{\Der}[2]{\frac{\text{d}#1}{\text{d}#2}}

\newcommand{\mathd}{\mathrm{d}}

\newcommand{\BF}[1]{\boldsymbol{#1}}

\newtheorem*{remark}{\bf Remark}

\definecolor{UoM}{RGB}{100,30,160}
\definecolor{newgreen}{RGB}{0,150,0}

\newcommand{\RED}{} 

\begin{document}

\title{Diffraction of acoustic waves by multiple semi-infinite arrays$^1$}
\author{M. A. Nethercote$^2$,~A. V. Kisil$^3$~and~R. C. Assier$^3$}
\footnotetext[1]{Accepted by The Journal of the Acoustical Society in America on 18/08/2023.}
\footnotetext[2]{Department of Applied Mathematics and Theoretical Physics, University of Cambridge, United Kingdom,\newline
	(mn598@cam.ac.uk.)}
\footnotetext[3]{Department of Mathematics, Alan Turing Building, University of Manchester, United Kingdom,\newline
	(anastasia.kisil@manchester.ac.uk, raphael.assier@manchester.ac.uk).}

\date{} 
\maketitle

\begin{abstract}
Analytical methods are fundamental in studying acoustics problems. One of the important tools is the Wiener-Hopf method, which can be used to solve many canonical problems with sharp transitions in boundary conditions on a plane/plate. However, there are some strict limitations to its use, usually the boundary conditions need to be imposed on parallel lines (after a suitable mapping). Such mappings exist for wedges with continuous boundaries, but for discrete boundaries, they have not yet been constructed. In our previous article, we have overcome this limitation and studied the diffraction of acoustic waves by a wedge consisting of point scatterers. Here, the problem is generalised to an arbitrary number of periodic semi-infinite arrays with arbitrary orientations. This is done by constructing several coupled systems of equations (one for every semi-infinite array) which are treated independently. The derived systems of equations are solved using the discrete Wiener–Hopf technique and the resulting matrix equation is inverted using elementary matrix arithmetic. Of course, numerically this matrix needs to be truncated, but we are able to do so such that thousands of scatterers on every array are included in the numerical results. Comparisons with other numerical methods are considered, and their strengths/weaknesses are highlighted.
\end{abstract}

\section{Introduction}
Analytical solutions for acoustic/electromagnetic wave scattering problems by different combinations of finite, infinite and semi-infinite plates/arrays \cite{LawrieAbrahams2007} are of special interest. These are difficult problems, and we will briefly review some of the work on this subject. Generalising the famous Sommerfeld’s half-plane problem (consisting of one semi-infinite plate), \citet{Heins1948part1,Heins1948part2} uses the Wiener--Hopf (WH) technique to find the exact solution for the problem where an electric-polarised wave is incident on a pair of parallel semi-infinite plates, symbolising a receiving and transmitting antenna for parts I and II respectively. This WH formulation has been extended to a matrix version in \cite{Abrahams1988,Abrahams1988part2,Abrahams1990}, for the equivalent acoustic problem with a pair of staggered plates, and in \cite{DSJones1986} for three semi-infinite planes. Also, by employing appropriate mappings/transformations, it was possible to use the matrix WH techniques for wedges \cite{AVShanin1998,DanieleZich2014,WedgeReview}. Additionally, since finding exact solutions is difficult for the matrix WH technique, there have been many developments on approximate and asymptotic factorisation techniques \cite{RogosinMishuris2016,KisilAyton2018,Kisil2018}. Alternatively, researchers resort to numerical and asymptotic schemes in their models (see \cite{PeakeCooper2001,Kirby2008,AdamsCrasterGuenneau2008,CrasterGuenneauAdams2009} for waveguides and ducts for example). 

One of the advantages of the WH technique is that a solution for a scattering problem can also be used in aeroacoustics setting to study the interaction of plates with gusts. In particular, \cite{Peake1992} considered an infinite staggered cascade of finite thin blades which were aligned with a uniform subsonic mean flow, and found an iterative solution by an infinite sequence of coupled, semi-infinite WH problems in the high frequency limit. More recently, this work has been extended for non-aligned mean flow \cite{PeakeKerschen1997,PeakeKerschen2004}, thin aerofoils \cite{Baddoo2018,Baddoo2020} and for high staggering angles \cite{MaierhoferPeake2020,MaierhoferPeake2022}.

In an elastic setting, the scattering and localisation of flexural waves on an elastic Kirchhoff plate is another important problem, especially determining waves that are blocked or trapped \cite{Movchan2009,Haslinger2014,Haslinger2016}. \citet{JonesMovchan2017} considered a pair of parallel semi-infinite gratings of rigid pins on a Kirchhoff plate. A follow-up article \cite{HaslingerJonesMovchan2018} studied problems with configurations of parallel semi-infinite gratings. In particular, it featured a problem where four parallel semi-infinite gratings form a waveguide in a herringbone pattern. Both these articles form and use the solution to the discrete WH functional equation but also identify trapped modes from analysing the kernel. The latter of the two articles also assumes that for each side of the waveguide, the two gratings are closely spaced which allowed them to use a dipole approximation.

These diffraction problems can also be considered on a lattice governed by a discrete Helmholtz equation. For instance, waves diffracted in a square lattice by a semi-infinite crack or a semi-infinite rigid constant was solved recently using discrete WH technique in \cite{BLSharma_SIAM2015a} and \cite{BLSharma_WM2015} respectively. These two problems are analogous to the classic Sommerfeld’s half-plane problem (sound-hard and sound-soft respectively). This work has been extended to two staggered semi-infinite cracks \cite{MauryaSharma2019} which was not solved exactly but asymptotically due to the notorious difficulties in factorising a matrix WH kernel.

In this article, we are interested in arbitrary arrays of small Dirichlet cylinders within a continuum. This means that we use the continuous Helmholtz equation subject to boundary conditions imposed on a discrete set of scatterers. The semi-infinite array (the analogue of Sommerfeld’s half-plane) problem was solved using the discrete WH technique long ago \cite{HillsKarp1965,LintonMartin2004}.  But there has been little work to generalise this to other configurations of arrays as was the case in the continuous and the discrete case outlined above. In our previous work \cite{HWpaperI}, we combined two semi-infinite arrays to form a wedge. Unlike the continuous boundary wedge \cite{WedgeReview}, this did not lead to a matrix WH problem due to the difficulties in finding an appropriate mapping. Instead, we considered two coupled systems of equations which were solved using the discrete WH technique, followed by an effective numerical iterative procedure. In this article, we will study problems involving any number of independent semi-infinite arrays comprised of equidistant scatterers. This is very general since the position and orientation of the arrays are arbitrary, which allows us to model many types of interesting problems. As in \cite{HWpaperI}, we will use the WH technique for each of the arrays and then couple them together. But this time, the coupling is encoded directly in the matrix inversion which means that there is no need for an iterative scheme.

For any arbitrary configuration of scatterers, one can use numerical techniques to determine the scattering behaviour. Some examples of these techniques include finite element methods, a T-matrix reduced order model \cite{HawkinsGanesh2017,tmatsolver} and a least square collocation approach that was used in \cite{ChapmanHewettTrefethen2015} and \cite{HewettHewitt2016} to study the electrostatic and electromagnetic shielding by Faraday cages. While these methods are very efficient at modelling the interactions between individual scatterers, they do not work very well at modelling the infinite nature of periodic arrays.

The structure of the paper is as follows. We start by setting up and solving the Wiener-Hopf problem in section \ref{sec:WHT}, which results in a matrix equation that is inverted in section \ref{sec:matrix} to find the scattering coefficients. In section \ref{sec:two_arrays}, we proceed by looking into the special case with two semi-infinite arrays and link with the point scatterer wedge from \cite{HWpaperI}. We also analyse the determinants and condition numbers of the matrices involved in section \ref{sec:uniqueness} as well as discuss the use of fast multipole methods for efficient computations of their components in section \ref{sec:FMM}. Finally, we showcase several different test cases in section \ref{sec:test_cases} and compare the Wiener-Hopf solution with other numerical techniques in section \ref{sec:comparison} and highlight their strengths and weaknesses.

\section{Multiple semi-infinite arrays}
Viewed as a three-dimensional problem, the scatterers are all cylinders of infinite height, have a small radius and satisfy homogeneous Dirichlet boundary conditions. This problem can naturally be reduced to two dimensions for non-skew incidence and this is what we will be considering here. Throughout this article, we will exploit the methodology \RED{described} in \cite{HWpaperI}, since the point scatterer wedge \RED{(or line scatterer wedge when viewed in three dimensions)} can be considered as a particular case of what is presented here. 

Similarly to \cite{HWpaperI}, we are looking for time-harmonic solutions to the linear wave equation by assuming and then suppressing the time factor $e^{-i\omega t}$, where $\omega$ is the angular frequency, and use a polar coordinate system $(r,\theta)$ with the position vector given by $\BF{r}$. We let $\Phi$ be the total wave field and decompose it into an incident wave field $\Phi_{\textrm{I}}$ and the resulting scattered field $\Phi_{\textrm{S}}$ by the equation $\Phi=\Phi_{\textrm{I}}+\Phi_{\textrm{S}}$. \RED{Each} of these fields satisfy the Helmholtz equation with wavenumber $k$. The incident wave field $\Phi_{\textrm{I}}(\BF{r})$ takes the form of a unit amplitude plane wave given by,
\begin{align}\label{MSIA-PhiI}
\Phi_{\textrm{I}}=e^{-ikr\cos(\theta-\theta_{\textrm{I}})}
\end{align}
where $\theta_{\textrm{I}}$ is the incoming incident angle. An important assumption of this study is the use of Foldy's approximation \cite{Foldy1945,pmartin2006}: where we assume that the cylinders are isotropic point scatterers. That requires them to be small in comparison to the wavelength (i.e.\ $ka\ll1$) and this allows us to write the scattered field $\Phi_{\textrm{S}}$ in the form of a monopole expansion. This article is focused on the problem of an incident wave scattered by $\mathcal{J}$ arbitrary periodic semi-infinite arrays. The $j^{\textrm{th}}$ array starts at an arbitrary position $\BF{R}^{(j)}_0$, makes an arbitrary angle $\alpha_j$ with the $x$-axis, and the scatterers have radius $a_j>0$ and are arranged with a spacing $s_j>0$. The position of the $n^{\textrm{th}}$ scatterer in the $j^{\textrm{th}}$ array is hence given by,
\begin{align}\nonumber\label{MSIA-R_n}
\BF{R}^{(j)}_n&=\BF{R}^{(j)}_0+ns_j(\cos(\alpha_j)\BF{\hat{x}}+\sin(\alpha_j)\BF{\hat{y}}),\quad n=0,1,2,...\\
\BF{R}^{(j)}_0&=R^{(j)}_0(\cos(\theta^{(j)}_0)\BF{\hat{x}}+\sin(\theta^{(j)}_0)\BF{\hat{y}})
\end{align}
where $(\BF{\hat{x}},\BF{\hat{y}})$ are the unit basis vectors of a Cartesian coordinate system. \RED{Note that $a_j$, $s_j$ and $R^{(j)}_0$ are all lengths that are measured in metres [m] and $k$ is in [m$^{-1}$]. These units will be omitted from now on.}
\begin{figure}[h]\centering
\includegraphics[width=0.35\textwidth]{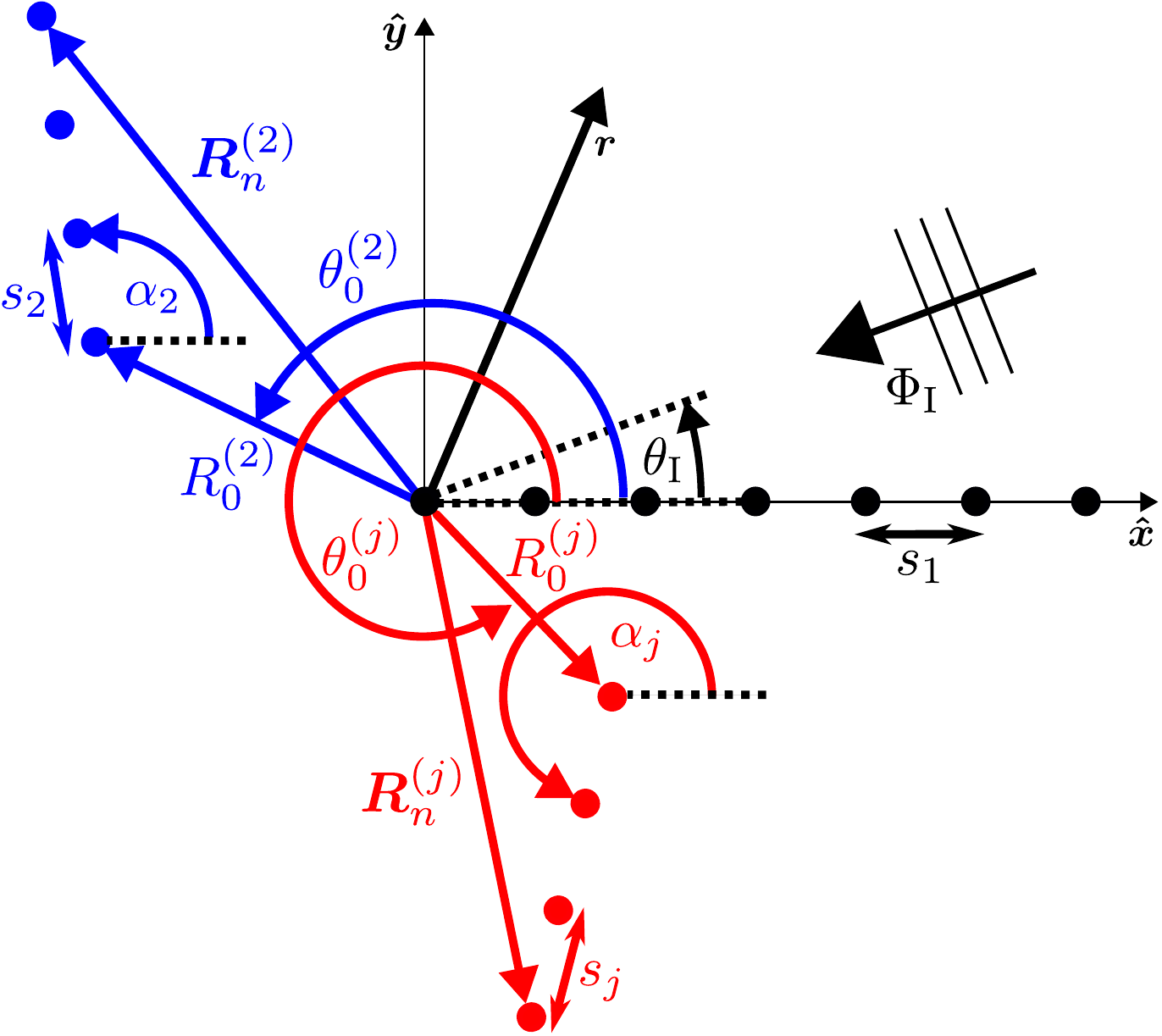}
\caption{Diagram of a plane wave interacting with multiple arbitrary semi-infinite arrays. For simplicity here, the first array is positioned on the positive $x$-axis (i.e.\ $R^{(1)}_0=\alpha_1=0$).}
\label{fig:MSIA-diagram}
\end{figure}

We further introduce 
\RED{\begin{align}
\Lambda^{(j,\ell)}(m,n)=|\BF{R}^{(j)}_m-\BF{R}^{(\ell)}_n|
\end{align}}
as the distance between the $m^{\textrm{th}}$ scatterer on the $j^{\textrm{th}}$ array and the $n^{\textrm{th}}$ on the $\ell^{\textrm{th}}$ array. 
This distance function satisfies the identity $\Lambda^{(j,\ell)}(m,n)=\Lambda^{(\ell,j)}(n,m)$. It is important to note that while we allow the arrays to cross, we do \emph{not} want the scatterers to overlap. This means that we need the condition $a_j<s_j/2$ for all $j=1,2,...\mathcal{J}$ to prevent overlapping between scatterers belonging to the $j^{\textrm{th}}$ array and $a_j+a_\ell<\Lambda^{(j,\ell)}(m,n)$ for all $m,n\in\ZZ$ and $j,\ell=1,2,...\mathcal{J}$ which prevents overlapping between the $j^{\textrm{th}}$ and $\ell^{\textrm{th}}$ arrays.

Using Foldy's approximation, the scattered field $\Phi_{\textrm{S}}$ is written in the form of a monopole expansion. This means that the total field $\Phi$ is given by
\begin{align}\label{MSIA-gensol}
\Phi(\BF{r})=\Phi_{\text{I}}+\sum_{j=1}^\mathcal{J}\sum_{n=0}^\infty A^{(j)}_nH^{(1)}_0(k|\BF{r}-\BF{R}^{(j)}_n|),
\end{align}
where $A^{(j)}_n$ is the scattering coefficient associated with the $n^{\textrm{th}}$ scatterer of the $j^{\textrm{th}}$ array. To obtain the systems of equations, we use the procedure described in \cite[eqns (3.4), (3.6)]{HWpaperI} for the point scatterer wedge. As a result, we obtain $\mathcal{J}$ systems of infinitely many equations governing the scattering coefficients. The $m^{\textrm{th}}$ equation ($m\geq0$) of the $j^{\textrm{th}}$ system ($j\in\{1,2,...\mathcal{J}\}$) is given by 
\begin{align}\label{MSIA-A-system}
&A^{(j)}_mH^{(1)}_0(ka_j)+\sum_{\substack{n=0\\n\neq m}}^{\infty}A^{(j)}_nH^{(1)}_0(ks_j|m-n|)
=-e^{i\BF{k}\cdot\BF{R}^{(j)}_m}-\sum_{\substack{\ell=1\\\ell\neq j}}^\mathcal{J}\sum_{n=0}^{\infty}A^{(\ell)}_nH^{(1)}_0\left(k\Lambda^{(j,\ell)}(m,n)\right)
\end{align}
where $\BF{k}\cdot\BF{R}^{(j)}_m=-kR^{(j)}_0\cos(\theta^{(j)}_0-\theta_{\text{I}})-ks_jm\cos(\alpha_j-\theta_{\text{I}})$ and $\BF{k}=-k(\cos(\theta_{\text{I}})\BF{\hat{x}}+\sin(\theta_{\text{I}})\BF{\hat{y}})$ is the incident wavevector. \RED{It is important to note that due to Foldy's approximation, the boundary conditions are only approximately satisfied. Regarding energy, \cite[Chp.~8]{pmartin2006} showed that it is conserved provided that
\begin{align}\label{MSIA-CoE}
\left|g^{(j)}_n\right|^2+\RE{g^{(j)}_n}=0,
\end{align}
where $g^{(j)}_n=A^{(j)}_n/\Phi^{(j)}_n(\BF{R}^{(j)}_n)$ with $$\Phi^{(j)}_n(\BF{r})=\Phi(\BF{r})-A^{(j)}_nH^{(1)}_0(k|\BF{r}-\BF{R}^{(j)}_n|).$$
Following \cite{HillsKarp1965}, for $ka_j\ll1$, we approximate $\Phi^{(j)}_n(\BF{R}^{(j)}_n)$ by $\Phi^{(j)}_n$ on the boundary of the associated scatterer to obtain $g^{(j)}_n=-(H^{(1)}_0(ka_j))^{-1}$. This leads to a discrepancy in \eqref{MSIA-CoE} of the order $O\Big(\left(\tfrac{ka_j}{\ln(ka_j)}\right)^2\Big)$, which is very small when $ka_j\ll1$.  Though we do not do this in the present work, it is possible for \eqref{MSIA-CoE} to be satisfied exactly by considering only the leading order approximation of $g^{(j)}_n$ as did \citet{LintonMartin2004}, but this does not change the solution to leading order.
}

\subsection{Solving the $j^{\textrm{th}}$ system of equations}\label{sec:WHT}
To solve the system of equations \eqref{MSIA-A-system} for a specific $j$, we use the discrete analogue of the WH technique. We start by extending \eqref{MSIA-A-system} for negative $m$ using some unknown coefficients $F_{j,m}$ and state that $A^{(j)}_m=0$ for $m<0$, 
\begin{align}\label{MSIA-begin-system}
&A^{(j)}_mH^{(1)}_0(ka_j)+\sum_{\substack{n=0\\n\neq m}}^{\infty}A^{(j)}_nH^{(1)}_0(ks_j|m-n|)
\!=\!\begin{dcases}\!-e^{i\BF{k}\cdot\BF{R}^{(j)}_m}\!-\!\sum_{\substack{\ell=1\\\ell\neq j}}^\mathcal{J}\!\sum_{n=0}^{\infty}\!A^{(\ell)}_nH^{(1)}_0\!\left(\!k\Lambda^{(j,\ell)}(m,n)\!\right) ,&m\geq0,\\
F_{j,m},&m<0.\end{dcases}
\end{align}
Here, the $A^{(j)}_m$ scattering coefficients are the unknowns to find and all others are assumed to be known. Noting that the forward \RED{and inverse Z--transform for any sequence $G_m$ are given by,
\begin{align}\label{MSIA-Z-transform}
G(z)=\!\sum_{m=-\infty}^\infty\! G_mz^m,\ G_m=\!\frac{1}{2\pi i}\oint_{C}\! G(z)z^{-m-1}\text{d}z,
\end{align}}
we apply it to \eqref{MSIA-begin-system} to obtain the Wiener--Hopf equation
\begin{align}\label{MSIA-WHE}
K_j(z)A_j^+(z)=&F_{j,\text{pole}}^+(z)+F_j^-(z)+\sum_{\substack{\ell=1\\\ell\neq j}}^{\mathcal{J}}F_{\ell,A}^+(z),
\end{align}
where $A_j^+(z)$ is the Z--transform of the unknown scattering coefficients of the $j^{\textrm{th}}$ array;
\begin{align}\label{MSIA-A+def}
A_j^+(z)=\sum_{m=-\infty}^{\infty} A^{(j)}_mz^m=\sum_{m=0}^{\infty} A^{(j)}_mz^m.
\end{align}
As in \cite{HWpaperI}, it is useful to assume that $k$ has a small positive imaginary part to help with the convergence of the Z--transform. We also define the two regions, 
\begin{align}\nonumber\label{MSIA-regions}
\Omega_j^+&=\left\{z\in\CC: |z| < e^{-\IM{k}s_j\cos(\alpha_j-\theta_{\textrm{I}})}\right\},\\ 
\Omega_j^-&=\left\{z\in\CC: |z| > e^{-\IM{k}s_j}\right\},
\end{align}
in which a function with a $+$ or $-$ superscript is analytic. In these Wiener--Hopf problems, a crucial function is the Wiener--Hopf \textit{kernel} $K_j(z)$ given by
\begin{align}\label{MSIA-Kj(z)}
K_j(z)=H^{(1)}_0(ka_j)+\sum_{\ell=1}^\infty (z^\ell+z^{-\ell})H^{(1)}_0(ks_j\ell),
\end{align}
for $j=1,2,...\mathcal{J}$, which has the exact same definition and properties as in \cite[eq (2.15)]{HWpaperI}, including the important identity $K_j(z)=K_j(1/z)$ and the singular points $z=e^{\pm iks_j}$. Furthermore, the kernel is analytic and zero-free on an annulus which contains $\Omega_j^+\cap\Omega_j^-$. Since $K_j(z)$ is a slow-convergent infinite series, it is very impractical for numerical evaluation. To counter this, there are alternative methods of evaluation, including the use of the method of tail-end asymptotics \cite{Lynott_19} and rewriting the Schl\"{o}milch series to a fast-convergent version \cite{Linton1998,Linton2006,Linton2010} (see also the appendix in \cite{HWpaperI} for specifics).

The three forcing terms on the right-hand side of \eqref{MSIA-WHE} are defined by
\begin{align}\label{MSIA-F(z)}
\nonumber F_{j,\text{pole}}^+(z)&\!=\!\frac{e^{i\BF{k}\cdot\BF{R}^{(j)}_0}}{ze^{-iks_j\cos(\alpha_j-\theta_{\text{I}})}-1},\\
\nonumber F_{\ell,A}^+(z)&\!=\!-\!\sum_{m=0}^\infty\sum_{n=0}^{\infty}\!A^{(\ell)}_{n}z^mH^{(1)}_0\!\left(\!k\Lambda^{(j,\ell)}(m,n)\!\right)\!,\\
F_j^-(z)&\!=\!\sum_{m=-\infty}^{-1}F_{j,m}z^m.
\end{align}
Note that by design, $F_j^-(z)=O\left(\frac{1}{z}\right)$ as $|z|\rightarrow\infty$. 

\begin{figure}[ht]\centering
\includegraphics[width=0.4\textwidth]{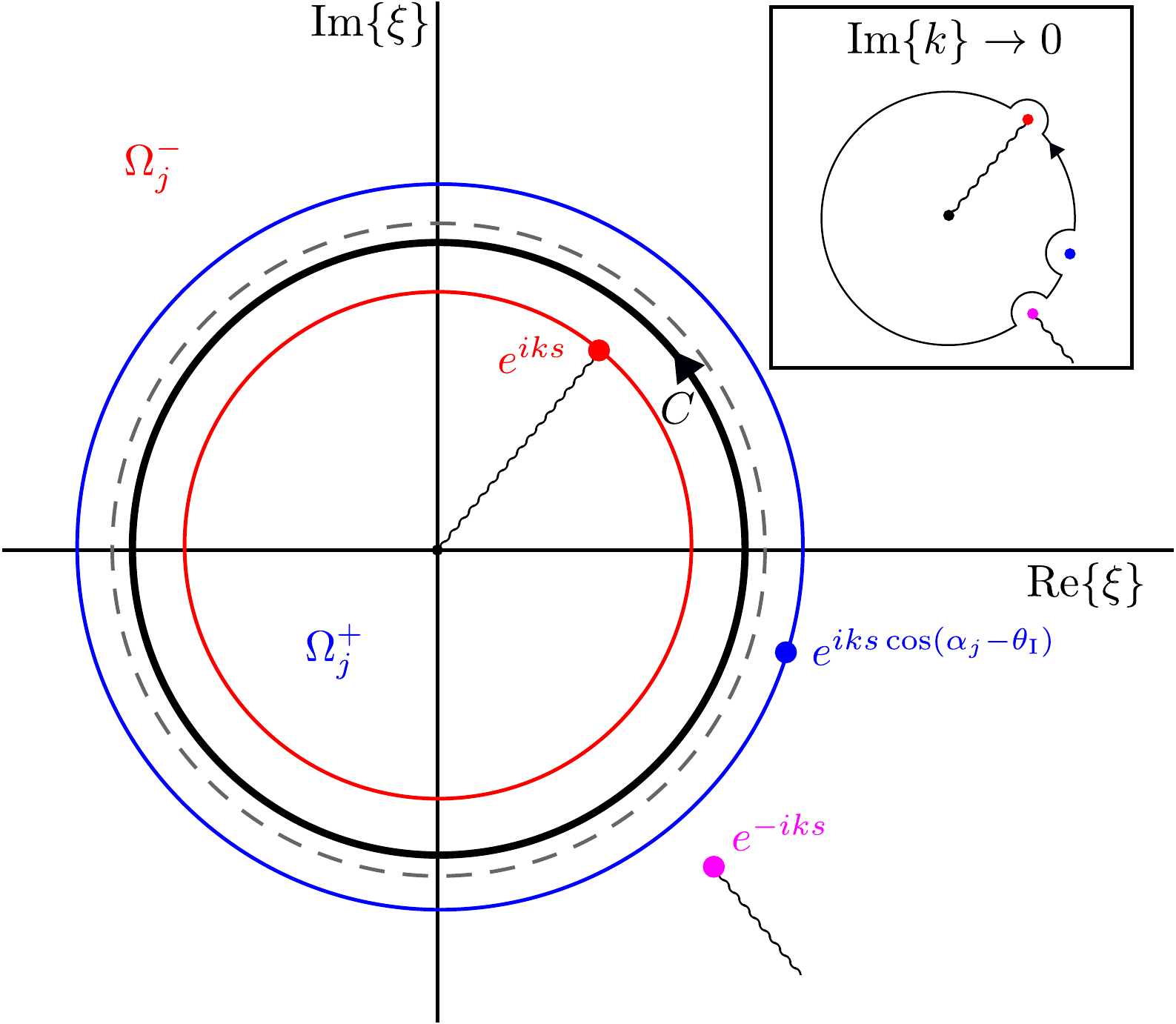}
\caption{Diagram of the integration contour $C$ on the $\xi$ complex plane. Here, the border of the regions $\Omega_j^\pm$ are shown as blue and red circles respectively and the grey dashed circle is the unit circle $|\xi|=1$. The smaller diagram is the limiting case when the imaginary part of $k$ tends to zero.}
\label{fig:C-contour}
\end{figure}

We proceed with the WH technique by factorising $K_j(z)$ in the exact same way as in \cite[eq (2.15)]{HWpaperI}. This means writing $K_j(z)=K_j^+(z)K_j^-(z)$ in such a way that the two factors satisfy $K_j^-(1/z)=K_j^+(z)$ and are defined by Cauchy's integral formulae \eqref{MSIA-Kj^+(z)}-\eqref{MSIA-K0}. 
\begin{align}
\label{MSIA-Kj^+(z)}\ln(K_j^+(z))=\ln(K_{j}^{0})-&\frac{1}{2\pi i}\int_{C}\frac{\ln(K_j(\xi))}{\xi-1/z}\text{d}\xi,\\
\label{MSIA-Kj^-(z)}\ln(K_j^-(z))=\ln(K_{j}^{0})-&\frac{1}{2\pi i}\int_{C}\frac{\ln(K_j(\xi))}{\xi-z}\text{d}\xi,\\
\label{MSIA-K0}\ln(K_{j}^{0})=\ln(K_j^+(0))=&\frac{1}{4\pi i}\int_{C}\frac{\ln(K_j(\xi))}{\xi}\text{d}\xi.
\end{align}
Here, the integration contour $C$ (see FIG.\ \ref{fig:C-contour}) is the anticlockwise circular path contained inside $\Omega_j^+\cap\Omega_j^-$ on the $\xi$ complex plane. Additionally, $C$ will also run radially below the pole at $\xi=z^{\pm1}$. Both kernel factors $K_j^\pm(z)$ are also analytic and zero-free inside the regions $\Omega_j^\pm$.
Now let us divide \eqref{MSIA-WHE} by $K_j^-(z)$ to obtain,
\begin{align}\label{MSIA-WHE-fac}
K_j^+(z)A_j^+(z)=&\frac{F_{j,\text{pole}}^+(z)}{K_j^-(z)}+\frac{F_j^-(z)}{K_j^-(z)}+\sum_{\substack{\ell=1\\\ell\neq j}}^{\mathcal{J}}\frac{F_{\ell,A}^+(z)}{K_j^-(z)}.
\end{align}
Next, we sum-split the pole term using the pole removal technique \RED{on the pole at $z_j:=e^{iks_j\cos(\alpha_j-\theta_{\textrm{I}})}$} to get,
\begin{align*}
\frac{F_{j,\text{pole}}^+(z)}{K_j^-(z)}\!=\!\frac{F_{j,\text{pole}}^+(z)}{K_j^-(\RED{z_j})}
\!+\!F_{j,\text{pole}}^+(z)\left[\frac{1}{K_j^-(z)}\!-\!\frac{1}{K_j^-(\RED{z_j})}\right]\!,
\end{align*}
where the first and second terms \RED{on the right hand side} are analytic in $\Omega_j^+$ and $\Omega_j^-$ respectively. We can also sum-split all of the $\tfrac{F_{\ell,A}^+(z)}{K_j^-(z)}$ terms for every $\ell$ by first noting that it can be rewritten as a Laurent series given by
\begin{align}\label{MSIA-FA-laurent}
\frac{F_{\ell,A}^+(z)}{K_j^-(z)}&=D_{j,\ell}(z)=\sum_{n=-\infty}^{\infty}D_{j,\ell,n}z^n,\ \text{where}\ D_{j,\ell,n}=\frac{1}{2\pi i}\int_{C}D_{j,\ell}(z)z^{-n-1}\mathd z,
\end{align}
because it is analytic on an annulus region of $z$ containing $\Omega_j^+\cap\Omega_j^-$. The sum-split of $D_{j,\ell}(z)$ is trivial\RED{:
\begin{align}\label{MSIA-split-D}
D_{j,\ell}^+(z)=\sum_{n=0}^{\infty}D_{j,\ell,n}z^n,\ D_{j,\ell}^-(z)=\sum_{n=-\infty}^{-1}D_{j,\ell,n}z^n.
\end{align}}
After the sum-split, we obtain the final WH equation,
\begin{align}
\label{MSIA-WHE-LHS}&K_j^+(z)A_j^+(z)\!-\!\frac{F_{j,\text{pole}}^+(z)}{K_j^-(\RED{z_j})}\!-\!\sum_{\substack{\ell=1\\\ell\neq j}}^{\mathcal{J}}\!D_{j,\ell}^+(z)\\
\label{MSIA-WHE-RHS}&\!=\!F_{j,\text{pole}}^+(z)\!\left[\!\frac{1}{K_j^-(z)}\!-\!\frac{1}{K_j^-(\RED{z_j})}\!\right]
\!+\!\sum_{\substack{\ell=1\\\ell\neq j}}^{\mathcal{J}}\!D_{j,\ell}^-(z)\!+\!\frac{F_j^-(z)}{K_j^-(z)},
\end{align}
where the left (right) hand side are analytic in $\Omega_j^+$ and $\Omega_j^-$ respectively. The two sides of the WH equation are used to construct an entire function $\Psi_j(z)$ defined by
\begin{align}
\Psi_j(z)=\begin{cases}
\eqref{MSIA-WHE-LHS} & z\in\Omega_j^+,\\
\eqref{MSIA-WHE-RHS} & z\in\Omega_j^-,\\
\eqref{MSIA-WHE-LHS}=\eqref{MSIA-WHE-RHS} & z\in\Omega_j^+\cap\Omega_j^-.\\
\end{cases}
\end{align}
Each term of \eqref{MSIA-WHE-RHS} is $O\left(\frac{1}{z}\right)$ as $|z|\rightarrow\infty$, which implies that $\Psi_j$ is bounded and tends to zero at infinity as well as entire. Therefore, Liouville's theorem ensures that both \eqref{MSIA-WHE-LHS} and \eqref{MSIA-WHE-RHS} are equivalently zero, and then we obtain the solution for $A_j^+(z)$,
\begin{align}\label{MSIA-WHE-sol}
A_j^+(z)=&\frac{F_{j,\text{pole}}^+(z)}{K_j^+(z)K_j^-(\RED{z_j})}+\sum_{\substack{\ell=1\\\ell\neq j}}^{\mathcal{J}}\frac{D_{j,\ell}^+(z)}{K_j^+(z)}.
\end{align}
To recover the scattering coefficients we use the inverse Z--transform \eqref{MSIA-Z-transform},
which gives us,
\begin{align}\label{MSIA-sol-int}
A^{(j)}_m=&\frac{\RED{z_j}e^{i\BF{k}\cdot\BF{R}^{(j)}_0}}{2\pi iK_j^-(\RED{z_j})}\oint_{C} \frac{z^{-m-1}}{K_j^+(z)(z-\RED{z_j})}\text{d}z
+\sum_{\substack{\ell=1\\\ell\neq j}}^{\mathcal{J}}\frac{1}{2\pi i}\oint_{C} \frac{D_{j,\ell}^+(z)}{K_j^+(z)}z^{-m-1}\text{d}z
\end{align}
Next, we let the imaginary part of $k$ tend to zero and then the integration contour $C$ in \eqref{MSIA-sol-int} is an indented anticlockwise unit circle (by passing any singularities on the unit circle radially below) with the pole $z=0$ being the only singularity inside. The smaller diagram in FIG.\ \ref{fig:C-contour} illustrates $C$ in \eqref{MSIA-sol-int}, but note that we do not have the branch point at $e^{iks}$ here. To evaluate these integrals, we need to recall the identity $K_j^-(z)=K_j^+(1/z)$ and note the expansion of $\left(K_j^+(z)\right)^{-1}$ given by 
\begin{align*}
\frac{1}{K_j^+(z)}\!=\!\sum_{n=0}^\infty\!\lambda_{j,n} z^n,\ \text{where}\ \lambda_{j,n}\!=\!\frac{1}{n!}\Der{^n}{z^n}\!\left[\!\frac{1}{K_j^+(z)}\!\right]_{z=0}.
\end{align*}
For the first integral of \eqref{MSIA-sol-int}, the evaluation is equivalent to the associated semi-infinite array problem (see \cite[eqn (2.22)]{HWpaperI}), \RED{the solution of which is $A^{(j)}_{0,m}$}, where the extra factor $e^{i\BF{k}\cdot\BF{R}^{(j)}_0}$ accounts for the off-centre start of the array, \RED{
\begin{align}\label{MSIA-integral-left}
A^{(j)}_{0,m}&=\frac{z_je^{i\BF{k}\cdot\BF{R}^{(j)}_0}}{2\pi iK_j^-(z_j)}\oint_{C} \frac{z^{-m-1}}{K_j^+(z)(z-z_j)}\text{d}z,
=-\frac{e^{i\BF{k}\cdot\BF{R}^{(j)}_0}}{K_j^-(z_j)}\sum_{n=0}^m\lambda_{j,n}z_j^{n-m}.
\end{align}}
Each remaining term in \eqref{MSIA-sol-int} adds the interaction from the $\ell^{\textrm{th}}$ array to the $j^{\textrm{th}}$ array and is evaluated in much the same way as in \cite[eqns (3.25-3.27)]{HWpaperI}:\RED{
\begin{align}\label{MSIA-integral-right}
\frac{1}{2\pi i}\oint_{C} \frac{D_{j,\ell}^+(z)}{K_j^+(z)}z^{-m-1}\text{d}z
&=\sum_{n=0}^m\lambda_{j,m-n}D_{j,\ell,n}.
\end{align}}
The coefficients $D_{j,l,n}$ are given by \eqref{MSIA-FA-laurent}:
\RED{\begin{align*}
D_{j,\ell,n}\!=\!-\!\sum_{q=0}^{\infty}\sum_{p=0}^{\infty}\sum_{l=0}^{\infty}
\tfrac{\lambda_{j,p}A^{(\ell)}_{q}H^{(1)}_0\!\left(\!k\Lambda^{(j,\ell)}(l,q)\!\right)\!}{2\pi i}\!\int_{C}\!z^{l-n-p-1}\!\mathd z,
\end{align*}}
where the integral is non-zero only when $l-n-p=0$ which implies that,
\begin{align*}
D_{j,\ell,n}=-\sum_{q=0}^{\infty}\sum_{p=0}^{\infty}\lambda_{j,p}A^{(\ell)}_{q}H^{(1)}_0\left(k\Lambda^{(j,\ell)}(p+n,q)\right),
\end{align*}
and then the scattering coefficients are equal to
\RED{\begin{align}\label{MSIA-A-sum-sol}
&A^{(j)}_m=\ A^{(j)}_{0,m}
-\!\sum_{\substack{\ell=1\\\ell\neq j}}^\mathcal{J}\sum_{q=0}^{\infty}\sum_{p=0}^{\infty}\sum_{n=0}^m \! \lambda_{j,m-n}\lambda_{j,p}A^{(\ell)}_{q}\!H^{(1)}_0\!\!\left(\!k\Lambda^{(j,\ell)}(p+n,q)\!\right)\!.
\end{align}}

\subsection{Writing and solving the Wiener--Hopf solution as a matrix equation}\label{sec:matrix}
We can write the Wiener--Hopf solution \eqref{MSIA-A-sum-sol} in the form of an infinite matrix equation,
\begin{align}\label{MSIA-matrix-system}
\BF{A}^{(j)}&=\BF{A}^{(j)}_0-\sum_{\substack{\ell=1\\\ell\neq j}}^\mathcal{J}\mathcal{M}^{(j,\ell)}\BF{A}^{(\ell)},
\end{align}
where $\BF{A}^{(j)}$ and $\BF{A}^{(j)}_0$ are infinite column vectors of scattering coefficients with entries $A^{(j)}_m$ and \RED{$A^{(j)}_{0,m}$}
respectively. The infinite matrices $\mathcal{M}^{(j,\ell)}$ have entries
\begin{align}\label{MSIA-matrix-entries}
\mathcal{M}^{(j,\ell)}_{mq}\!=\!\sum_{p=0}^{\infty}\sum_{n=0}^m\!\lambda_{j,m-n}\lambda_{j,p}H^{(1)}_0\!\left(\!k\Lambda^{(j,\ell)}(p\!+\!n,q)\!\right)\!,
\end{align}
\RED{for $m,q\geq0$.} Putting together all values of $j$ gives us a system of matrix equations which can also be written in block matrix form,
\begin{align}\label{MSIA-matrixmatrix-system}
\begin{pmatrix}
\mathcal{I}&\mathcal{M}^{(1,2)}&\dots&\mathcal{M}^{(1,\mathcal{J})}\\
\mathcal{M}^{(2,1)}&\mathcal{I}&\dots&\mathcal{M}^{(2,\mathcal{J})}\\
\vdots&\vdots&\ddots&\vdots\\
\mathcal{M}^{(\mathcal{J},1)}&\mathcal{M}^{(\mathcal{J},2)}&\dots&\mathcal{I}
\end{pmatrix}\!\!\begin{pmatrix}
\BF{A}^{(1)}\\\BF{A}^{(2)}\\\vdots\\\BF{A}^{(\mathcal{J})}
\end{pmatrix}\!=\!
\begin{pmatrix}
\BF{A}^{(1)}_0\\\BF{A}^{(2)}_0\\\vdots\\\BF{A}^{(\mathcal{J})}_0
\end{pmatrix},
\end{align}
where $\mathcal{I}$ is the identity matrix. In principle, \eqref{MSIA-matrixmatrix-system} can be inverted to get,
\begin{align}\label{MSIA-matrixmatrix-solution}
\begin{pmatrix}
\BF{A}^{(1)}\\\BF{A}^{(2)}\\\vdots\\\BF{A}^{(\mathcal{J})}
\end{pmatrix}\!=\!
\begin{pmatrix}
\mathcal{I}&\mathcal{M}^{(1,2)}&\dots&\mathcal{M}^{(1,\mathcal{J})}\\
\mathcal{M}^{(2,1)}&\mathcal{I}&\dots&\mathcal{M}^{(2,\mathcal{J})}\\
\vdots&\vdots&\ddots&\vdots\\
\mathcal{M}^{(\mathcal{J},1)}&\mathcal{M}^{(\mathcal{J},2)}&\dots&\mathcal{I}
\end{pmatrix}^{-1}\!\!
\begin{pmatrix}
\BF{A}^{(1)}_0\\\BF{A}^{(2)}_0\\\vdots\\\BF{A}^{(\mathcal{J})}_0
\end{pmatrix}.
\end{align}
Note that to evaluate these matrices and their inverses in practice, we will need to truncate the summations in the entries as well as the block matrices themselves. We will do this by ensuring that all blocks in the big matrix are of the same size.

\subsection{Two semi-infinite arrays}\label{sec:two_arrays}
In this section, we would like to analyse the form of the inverse matrix in \eqref{MSIA-matrixmatrix-solution} in terms of its blocks. This is difficult to do analytically, especially for larger $\mathcal{J}$. Note that \eqref{MSIA-matrixmatrix-solution} reduces to the standard solution to the semi-infinite array problem if $\mathcal{J}=1$. Let us say that we have just two semi-infinite arrays (i.e.\ $\mathcal{J}=2$), then we have the following matrix system,
\begin{align}\label{MSIA-2array-system}\nonumber
\BF{A}^{(1)}&=\BF{A}^{(1)}_0-\mathcal{M}^{(1,2)}\BF{A}^{(2)},\\
\BF{A}^{(2)}&=\BF{A}^{(2)}_0-\mathcal{M}^{(2,1)}\BF{A}^{(1)}.
\end{align}
Here, the entries of the first terms are still given by \eqref{MSIA-integral-left}, 
and the entries of the matrices are given by,
\begin{align}\label{MSIA-2array-matrix-entries}\nonumber
\mathcal{M}^{(1,2)}_{mq}&\!=\!\sum_{p=0}^{\infty}\sum_{n=0}^m\! \lambda_{1,m-n}\lambda_{1,p}H^{(1)}_0\!\left(\!k\Lambda^{(1,2)}(p\!+\!n,q)\right),\\
\mathcal{M}^{(2,1)}_{mq}&\!=\!\sum_{p=0}^{\infty}\sum_{n=0}^m\! \lambda_{2,m-n}\lambda_{2,p}H^{(1)}_0\!\left(\!k\Lambda^{(1,2)}(q,p\!+\!n)\!\right)\!.
\end{align}
We can solve the system \eqref{MSIA-2array-system} using simple matrix arithmetic,
\begin{align}\nonumber\label{MSIA-2array-exact}
\BF{A}^{(1)}&=\left(\mathcal{I}-\mathcal{M}^{(1,2)}\mathcal{M}^{(2,1)}\right)^{-1}\left(\BF{A}^{(1)}_0-\mathcal{M}^{(1,2)}\BF{A}^{(2)}_0\right),\\
\BF{A}^{(2)}&=\BF{A}^{(2)}_0-\mathcal{M}^{(2,1)}\BF{A}^{(1)}.
\end{align}
There is also an equivalent alternative solution to \eqref{MSIA-2array-exact} where all the $1$'s and $2$'s have switched places. Although in theory, it is possible to write such formulae for $\mathcal{J}>2$, it quickly becomes very convoluted.

\begin{remark}
It is fairly simple to match these results with the specific case of the point scatterer wedge studied in \cite{HWpaperI}. That wedge configuration is produced from these parameter choices: $a_1=a_2=a$, $s_1=s_2=R_0^{(2)}=s$, $R_0^{(1)}=0$, $\theta_0^{(2)}=\alpha_2=-\alpha_1=-\alpha$ which gives the distance function $\Lambda^{(1,2)}(m,n)= s\left(m^2+(n+1)^2-2m(n+1)\cos(2\alpha)\right)^{\frac{1}{2}}$. This leads to the two WH kernels $(K_1(z)$ and $K_2(z))$ as well as the resulting coefficients $(\lambda_{1,n}$ and $\lambda_{2,n})$ becoming identical. There are two key differences left. One is that the second array here has one extra scatterer after truncation, i.e.\ the vector $\BF{A}^{(2)}$ has one extra entry. However, this extra entry can be neglected and consequently, one must ignore the last column in $\mathcal{M}^{(1,2)}$ and the last row in $\mathcal{M}^{(2,1)}$ as well. The other difference is that we have revoked the need for the iterative scheme and have instead inverted the matrix equation \eqref{MSIA-2array-system} directly. Note that the inverted matrix can be expanded, $\left(\mathcal{I}-\mathcal{M}^{(1,2)}\mathcal{M}^{(2,1)}\right)^{-1}= \mathcal{I}+\mathcal{M}^{(1,2)}\mathcal{M}^{(2,1)}+\left(\mathcal{M}^{(1,2)}\mathcal{M}^{(2,1)}\right)^2+...$, which can then be used to recover every iteration in the iterative scheme (although this requires the spectral radius $\rho\left(\mathcal{M}^{(1,2)}\mathcal{M}^{(2,1)}\right)<1$ to converge).
\end{remark}

\subsection{Uniqueness of Solution}\label{sec:uniqueness}
\begin{figure}[ht]\centering
\includegraphics[width=0.35\textwidth]{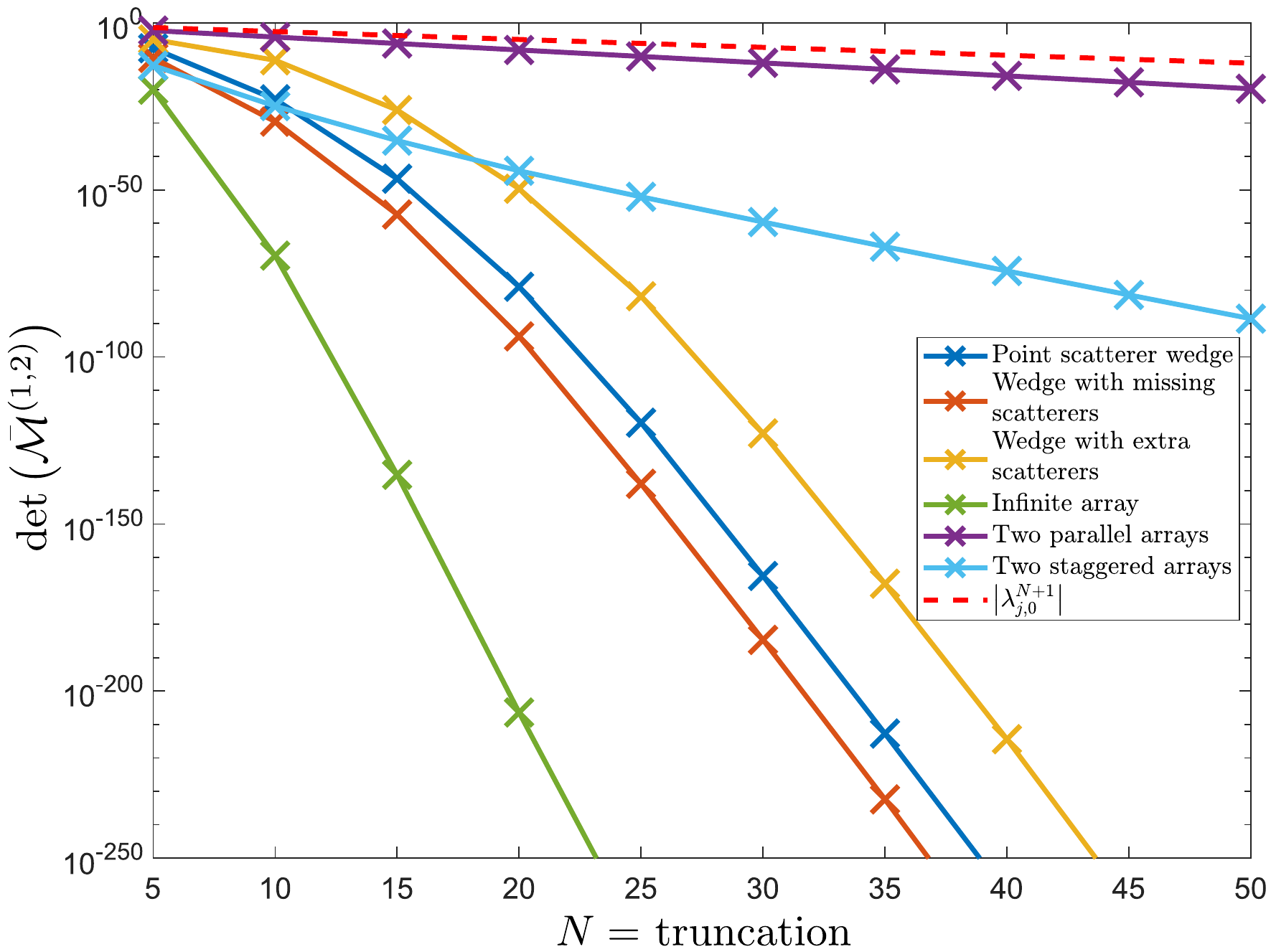}
\caption{Plot of the absolute value of $\det(\mathcal{M}^{(j,\ell)})$ w.r.t. the truncation $N$, compared with $\lambda_{j,0}^{N+1}$. For all cases, we have $k=5\pi$, $s_j=0.1$ and $a_j=0.001$ for all $j$ so the value of $\lambda_{j,0}$ is unchanged.}
\label{fig:MSIA-det}
\end{figure}

To be able to solve the matrix equation \eqref{MSIA-matrixmatrix-system}, we will need the determinant of that matrix to be non-zero. We have tried to find cases where the matrix is not invertible or show that it can always be inverted. While it seems to be the latter, it is clear from numerical experimentation that the individual block matrices are singular. The matrix entries can be written as follows, 
\begin{align}\nonumber\label{MSIA-barM}
\mathcal{M}^{(j,\ell)}&=\sum_{n=0}^m \lambda_{j,m-n}\bar{\mathcal{M}}^{(j,\ell)}_{nq},\\
\textrm{where}\ \ \bar{\mathcal{M}}^{(j,\ell)}_{nq}&=\sum_{p=0}^{\infty}\lambda_{j,p}H^{(1)}_0\left(k\Lambda^{(j,\ell)}(p+n,q)\right).
\end{align}
Say that we truncated the matrices at $N$ such that $\mathcal{M}^{(j,\ell)}$ is an $(N+1)$ by $(N+1)$ matrix, then by using Gaussian elimination, it can be shown that 
\begin{align}\label{MSIA-det-reduced}
\det(\mathcal{M}^{(j,\ell)})&=\lim_{N\rightarrow\infty}\left(\lambda_{j,0}^{N+1}\det(\bar{\mathcal{M}}^{(j,\ell)})\right),
\end{align}
For small $N$, the determinant of $\bar{\mathcal{M}}^{(j,\ell)}$ is not generally zero. From numerical experimentation however, we find that as $N$ increases, the extra eigenvalues (due to a bigger matrix) have very small absolute values and decay to zero. This means that $\det(\bar{\mathcal{M}}^{(j,\ell)})$ (being the product of all eigenvalues) decays to zero very fast (at least exponentially), which implies that $\det(\mathcal{M}^{(j,\ell)})$ decays to zero as well (assuming that $\lambda_{j,0}$ is sufficiently small in the worst case scenarios). In addition to this, the condition numbers for $\mathcal{M}^{(j,\ell)}$ should be infinite because these are singular matrices. FIG.\ \ref{fig:MSIA-det} plots the absolute value of $\det(\mathcal{M}^{(j,\ell)})$ w.r.t. $N$ for several different test cases. We find that the case with two parallel arrays has the slowest decay but it is still exponential.

The behaviour of the full matrix is naturally very different because the condition number for the matrix in \eqref{MSIA-matrixmatrix-system} is quite moderate (approximately $10^{0\text{-}2}$) and the determinants are non-zero for all cases that we tested for this article. If the system had a zero determinant, it could either be due to the modelling assumption being insufficient \cite{ArrayResonanceFoldyPaper} or due to a specific physical \RED{phenomenon}. Since a zero determinant implies a non-unique solution to the matrix equation, and considering that this mostly relies on the chosen geometry of the problem, one could conjecture that this could imply the presence of homogeneous (Rayleigh-Bloch) waves. However, since it was proven that semi-infinite arrays with Dirichlet boundary conditions cannot support these waves \cite{BonnetBenDhia1994,BonnetBenDhia2016}, we are not expecting a zero determinant.

\subsection{Computational optimisation}\label{sec:FMM}
\begin{figure}[ht]\centering
\includegraphics[width=0.75\textwidth]{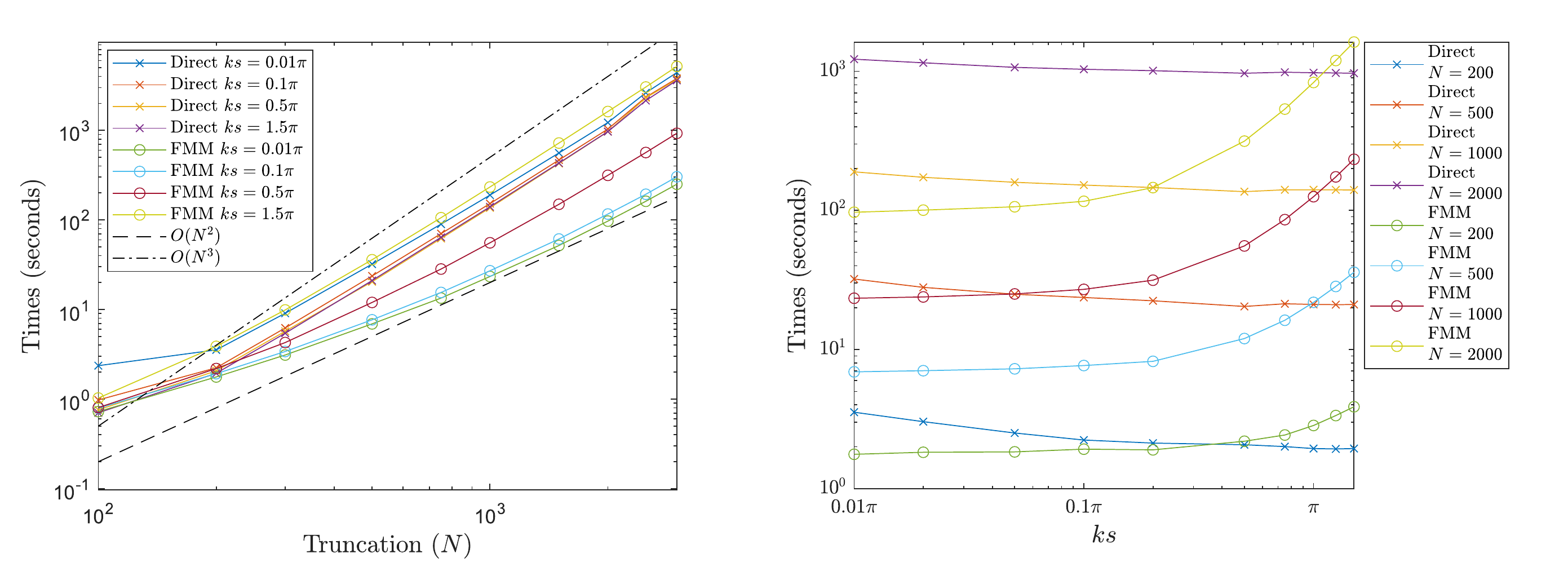}
\caption{Plots of the computation times to calculate the matrix in \eqref{MSIA-matrixmatrix-system} for the point scatterer wedge given in FIG.\ \ref{fig:test-cases}. On the left (resp. right) side, these times are plotted w.r.t. the truncation $N$ (resp. $ks$).}
\label{fig:MSIA-times}
\end{figure}

To calculate the matrices $\bar{\mathcal{M}}^{(j,\ell)}$, we can use the fast multipole methods (FMM) library which is accessible from \cite{FMM_library} (see also \cite{BeatsonGreengardFMM} for details on the algorithm). This method is very accurate at computing sums of the form \eqref{MSIA-barM}. For the matrices $\mathcal{M}^{(j,\ell)}$, the algorithm is able to reduce the computational cost with respect to the truncation ($N$) from $O(N^3)$ to as low as $O(N^2\ln(N))$. However, it is only able to do this if the values of $ks_j$ are sufficiently small for all $j$. To demonstrate this, FIG.\ \ref{fig:MSIA-times} plots the computation times to calculate the matrix in \eqref{MSIA-matrixmatrix-system} for the point scatterer wedge given in FIG.\ \ref{fig:test-cases}(a) \RED{(where $s_1=s_2=s$)}. This figure shows the difference between both methods w.r.t. the truncation $N$ and $ks$. On the left side, we see that for smaller $ks$ the computational order has been reduced. On the right side, we see that FMM becomes slower than using direct methods when $ks$ is approximately larger than $\pi$. Although there are developments for larger wavenumbers \cite{Greengardetal2006}, this has not been implemented in the library we used for the current work. As a result, we will only use FMM to calculate $\bar{\mathcal{M}}^{(j,\ell)}$ if the values of $ks_j$ are sufficiently small for all $j$. If not, then we will calculate it directly.

\begin{table}[ht]
\caption{The parameters for all semi-infinite arrays of the six test cases displayed in FIG.\ \ref{fig:test-cases}. \RED{In every case, we have $a_j=0.001$ for all arrays. The incident wave has the same parameters in all test cases, $k=5\pi$ and $\theta_{\textrm{I}}=\frac{\pi}{4}$ except for the last one (FIG.\ \ref{fig:test-cases}(f)) which has $k=7.5\pi$ instead.}}
\label{table:test-cases}
\begin{tabular}{c|c|c|c|c|c}
\hline\hline
{\footnotesize Case name} & $j$ & $s_j$ & $\alpha_j$ & $\theta^{(j)}_0$ & $R_0^{(j)}$ \\ \hline
{\footnotesize Point scatterer wedge} & $1$ & $0.1$ & $\frac{5\pi}{6}$ & $0$ & $0$\\ \cline{2-6}
{\footnotesize (FIG.\ \ref{fig:test-cases}(a))}& $2$ & $0.1$ & $-\frac{5\pi}{6}$ & $-\frac{5\pi}{6}$ & $0.1$\\ \hline
{\footnotesize Wedge with missing scatterers} & $1$ & $0.1$ & $\frac{5\pi}{6}$ & $\frac{5\pi}{6}$ & $0.3$\\ \cline{2-6}
{\footnotesize (FIG.\ \ref{fig:test-cases}(b))} & $2$ & $0.1$ & $-\frac{5\pi}{6}$ & $-\frac{5\pi}{6}$ & $0.3$\\ \hline
{\footnotesize Wedge with extra scatterers} & $1$ & $0.1$ & $\frac{5\pi}{6}$ & -$\frac{\pi}{6}$ & $0.45$\\ \cline{2-6}
{\footnotesize (FIG.\ \ref{fig:test-cases}(c))} & $2$ & $0.1$ & $-\frac{5\pi}{6}$ & $\frac{\pi}{6}$ & $0.45$\\ \hline
{\footnotesize Multilayered Faraday cage} & $1,2...7$ & $0.05$ & $\frac{(j-1)\pi}{6}$ & $\frac{(j-1)\pi}{6}$ & $0.1$\\ \cline{2-6}
{\footnotesize (FIG.\ \ref{fig:test-cases}(d))} & $8,9...12$ & $0.05$ & $\frac{(j-13)\pi}{6}$ & $\frac{(j-13)\pi}{6}$ & $0.1$\\ \hline
{\footnotesize Multiple infinite arrays} & $1,3...11$ & $0.1$ & $0$ & $-\frac{\pi}{2}$ & $0.1\left(\frac{j-1}{2}\right)$ \\ \cline{2-6}
{\footnotesize (FIG.\ \ref{fig:test-cases}(e) and (f))} & $2,4...12$ & $0.1$ & $\pi$ & $-\pi+\tan^{-1}\left(\frac{j}{2}-1\right)$ & \phantom{$\bigg|$}$0.1\sqrt{1+\left(\frac{j}{2}-1\right)^2}$\\
\hline\hline
\end{tabular}
\end{table}

\subsection{Test cases}\label{sec:test_cases}
\begin{figure}[t!]\centering
\includegraphics[width=0.8\textwidth]{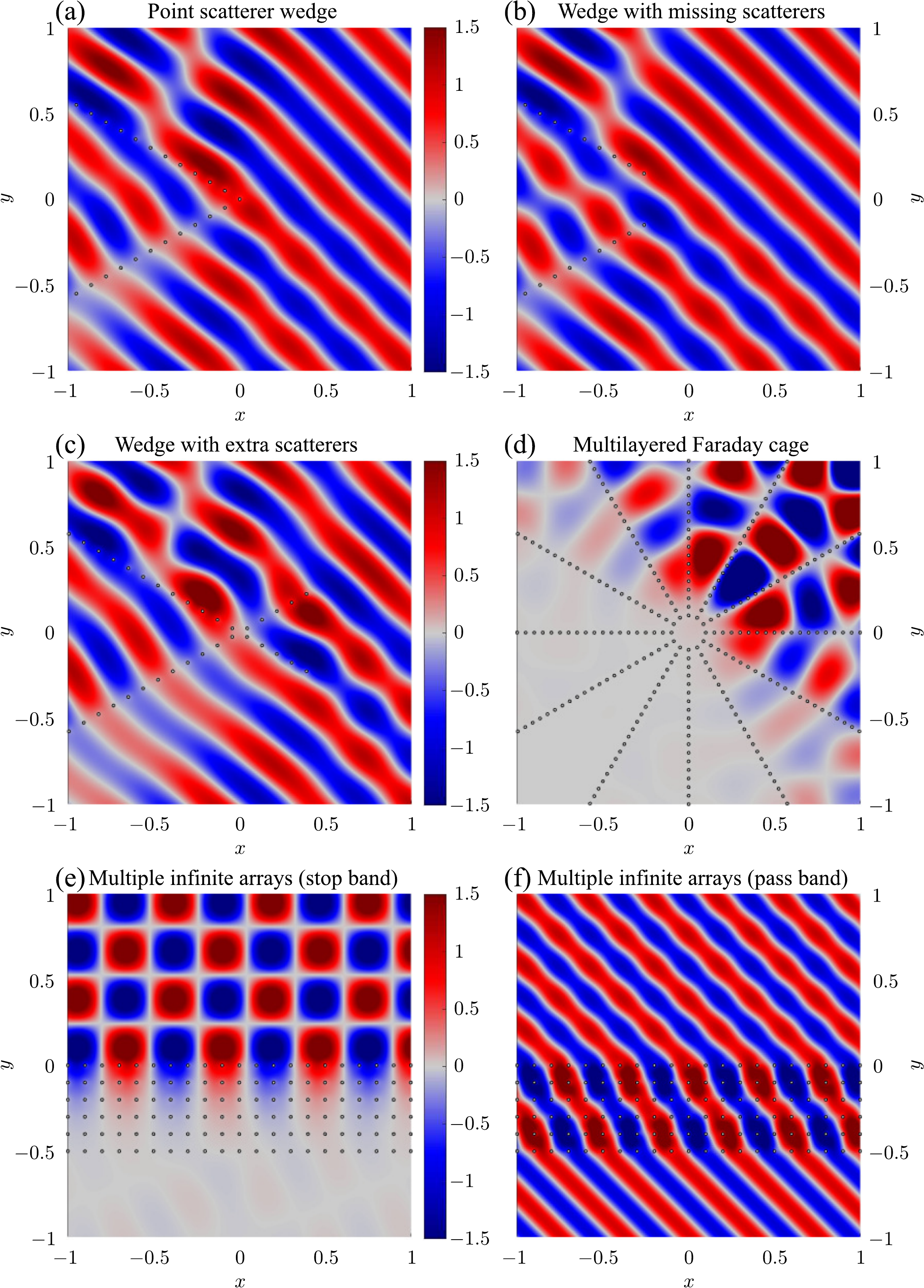}
\caption{Real part of total field for six different test cases. Here, the incident wave is given by the parameters $k=5\pi$ and $\theta_{\textrm{I}}=\frac{\pi}{4}$ (except for (f) plot which has $k=7.5\pi$), and the array parameters are given in Table \ref{table:test-cases}.}
\label{fig:test-cases}
\end{figure}

In this section, we consider and showcase several different examples of test cases in FIG.\ \ref{fig:test-cases} \RED{by plotting the real part of the total wave field.} We look at five different configurations of semi-infinite arrays, where the array parameters are given in Table \ref{table:test-cases}.

The first case (FIG.\ \ref{fig:test-cases}(a)) is a exact recreation of the same point scatterer wedge that was considered in Figure 6 of \cite{HWpaperI}. This time however, we are able to create gaps in the wedge interface (FIG.\ \ref{fig:test-cases}(b)) as well as add additional scatterers to one or both of the arrays (FIG.\ \ref{fig:test-cases}(c)), provided that we do not have overlapping scatterers. With this new generalised formulation, we can also consider cases with additional arrays. For example, we can have twelve outwardly pointing arrays where the ends are positioned to create a cage (FIG.\ \ref{fig:test-cases}(d)). This cage is able to shield the middle region from an incident wave with a low wavenumber. In particular, the sound pressure level inside the cage (given by the formula $20\log_{10}\left(\text{root mean square}(\Phi)\right)$) is approximately $-26.36$ compared to $0$ when there are no scatterers at all. This configuration is analogous to electrostatic and electromagnetic shielding problems using Faraday cages (see \cite{ChapmanHewettTrefethen2015,HewettHewitt2016}).

Although the WH technique is not typically used outside of semi-infinite problems, it is important to model the wave scattering by other types of arrays including circular arrays \cite{Martin2014}, infinite arrays with defects \cite{ThompsonLinton2008} and long finite arrays \cite{ThompsonLintonPorter2008}. Another case of special interest is determining the band-gap structure of doubly periodic lattices \cite{Botten2001}. For example, \citet{McIver2007} and \citet{KrynkinMcIver2009} study an infinite doubly periodic lattice of small scatterers with Dirichlet boundary conditions. With that in mind, the final case that we consider is a series of stacked infinite arrays to create a finitely thick doubly periodic lattice. The reason why we chose to look into this case is because we wanted to know if this configuration has similar properties to the fully periodic lattice as told by the band gap diagrams of \cite{KrynkinMcIver2009}. Specifically, we wanted to see the behaviour of an incident wave that can not penetrate the lattice (i.e.\ in a stop band) and the Bloch waves resulting from one that can (i.e.\ in a pass band). Choosing a wavenumber within a stop band (FIG.\ \ref{fig:test-cases}(e)) causes the incident wave to be almost fully reflected from the lattice. The alternative choice of a wavenumber within a pass band (FIG.\ \ref{fig:test-cases}(f)) does cause some reflection but most of the energy goes into the lattice and forms a Bloch wave inside before becoming a transmission out of the other side.

\subsection{Comparison with numerical methods}\label{sec:comparison}
Here, we seek to find a means of comparison with numerical methods that is both as fair and comprehensive as possible. We have three possible methods to compare with: finite element software COMSOL, a T-matrix solver (TMAT) by \cite{tmatsolver} and a least square collocation (LSC) method which was used in \cite{ChapmanHewettTrefethen2015,HewettHewitt2016} for solving Laplace's and Helmholtz's equation respectively. 

In previous articles \cite{LTpaper,HWpaperI}, we have extensively used COMSOL for comparison. However, there are limitations to the comparison since COMSOL is not able to find a solution with thousands of scatterers to fully compare with our method. Here, we will compare with results computed using TMAT or LSC where we are able to have the same number of scatterers and also the computation is performed on the same computer. We are also able to compute the scattering coefficients with these methods which is not possible with COMSOL.

The T-matrix software package provides an object-oriented implementation of an efficient reduced order model framework for modelling two- and three-dimensional wave scattering simulations. For the scattered field, TMAT uses the multipole expansion which is truncated depending on the scatterer size relative to the wavelength. TMAT also creates a Bessel function expansion of the incident field about every single scatterer which is truncated to the same number of terms as the multipole expansion. Following the construction of the T-matrix for every scatterer, a system of matrix equations is formed and solved for the scattering coefficients. However, the current version of this software restricts the truncation such that dipole coefficients must be included as well as monopole coefficients. 

The idea of the least squares collocation method is to create and solve an overdetermined matrix system to find the scattering coefficients of a truncated multipole expansion. Here, the known data is a collection of boundary data from a number of collocation points for every scatterer. We would need more collocation points per scatterer than the number of multipole terms to guarantee that the matrix system is overdetermined. The TMAT and LSC methods are both able to consider a large number of scatterers and despite the different methods, they are comparable in their results such that we only need to compare the WH method with one of them. We choose the LSC method because it will be a fairer comparison since we can restrict that one to monopole coefficients only.

If we were to compare each of the plots in FIG.\ \ref{fig:test-cases} with the equivalent determined by the LSC method, they would look identical without closer inspection. However, one can have a better idea of the differences between them by looking at the scattering coefficients produced by both methods. Let us consider a case with a single infinite array where the parameters are the same as the multiple infinite array case in Table \ref{table:test-cases}. The advantage here is that the infinite array problem has a known exact solution given by
\begin{align}\nonumber\label{MSIA-inf-array-sol}
A^{(1)}_m&=-\frac{e^{-iksm\cos(\theta_{\textrm{I}})}}{K(e^{iks\cos(\theta_{\textrm{I}})})},\quad m\geq0\\
A^{(2)}_m&=-\frac{e^{iks(m+1)\cos(\theta_{\textrm{I}})}}{K(e^{iks\cos(\theta_{\textrm{I}})})},\quad m\geq0
\end{align}
where $s=s_1=s_2$ and $K=K_1=K_2$. This means that we are able compare both the WH and LSC solutions with the exact solution as well as each other. 

\begin{figure}[ht]\centering
\includegraphics[width=0.49\textwidth]{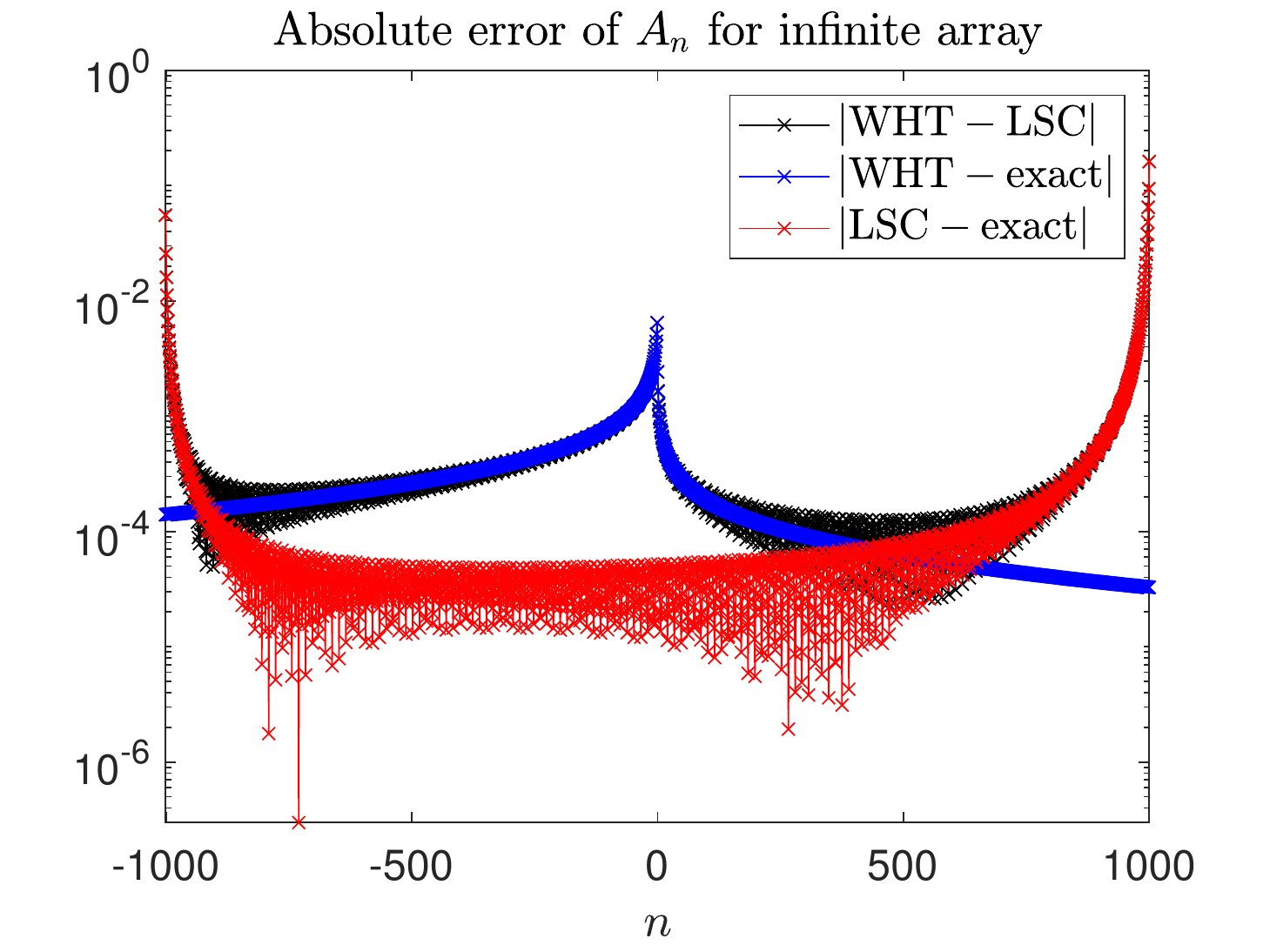}
\includegraphics[width=0.49\textwidth]{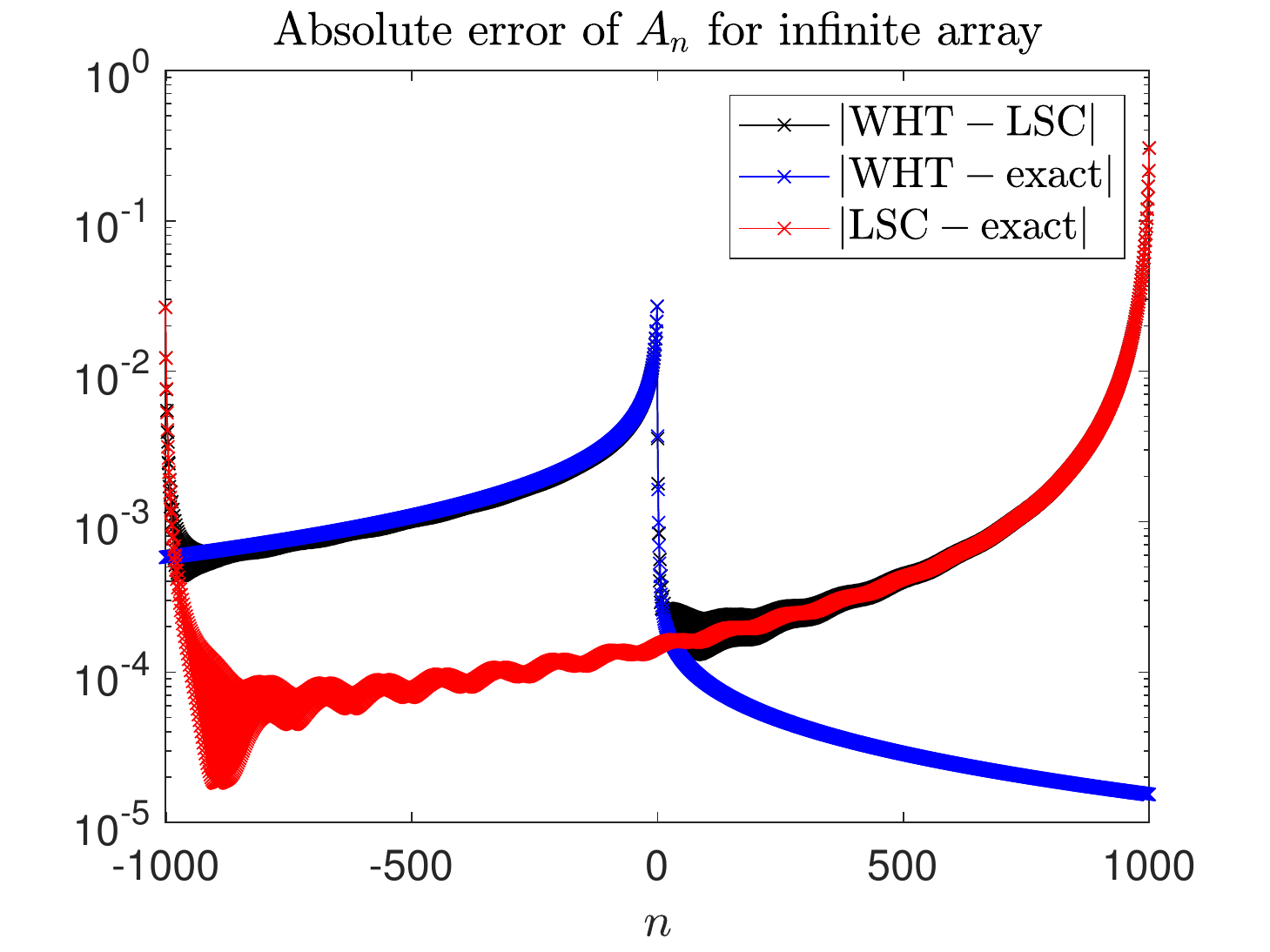}\\
\includegraphics[width=0.49\textwidth]{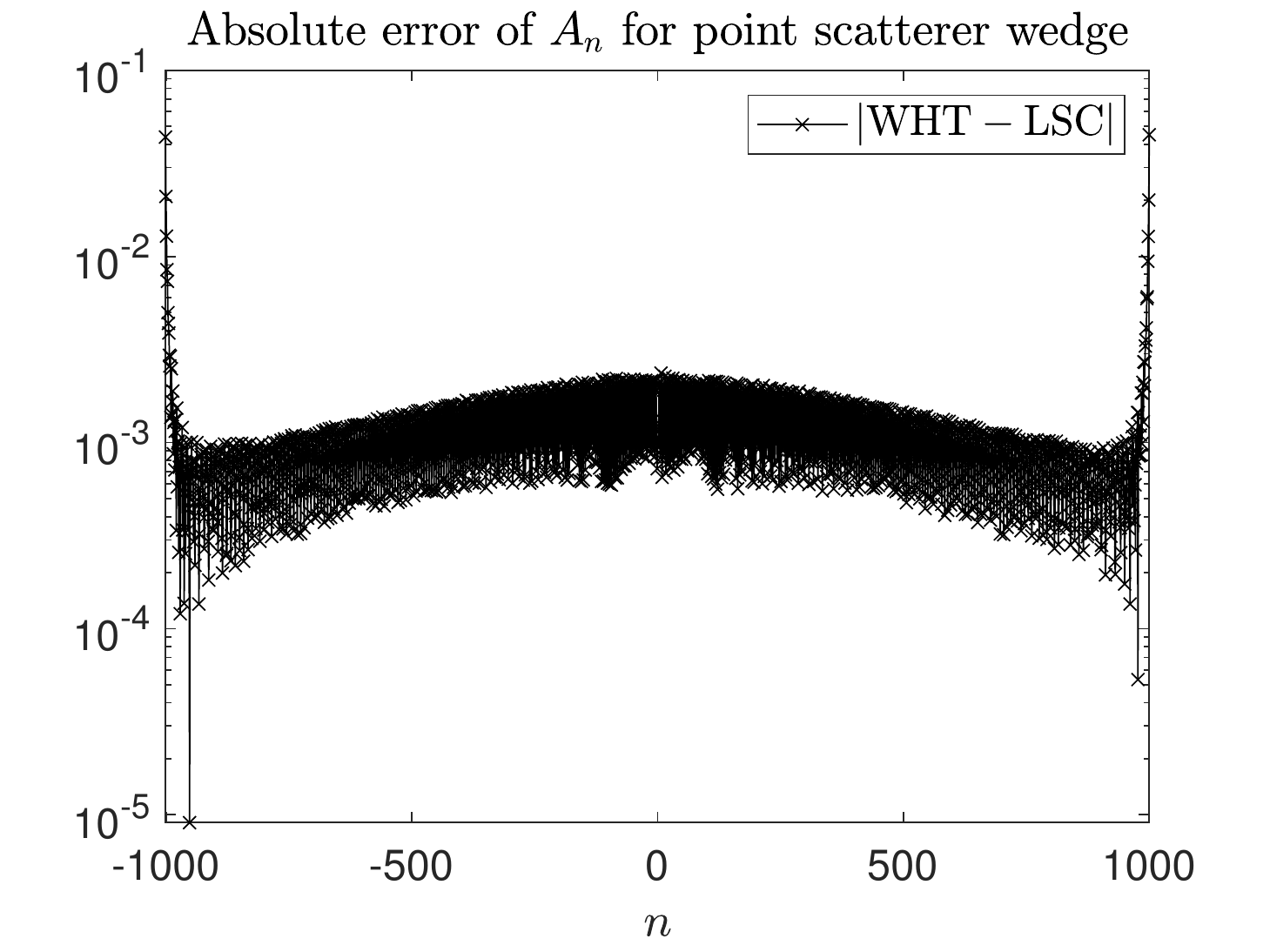}
\includegraphics[width=0.49\textwidth]{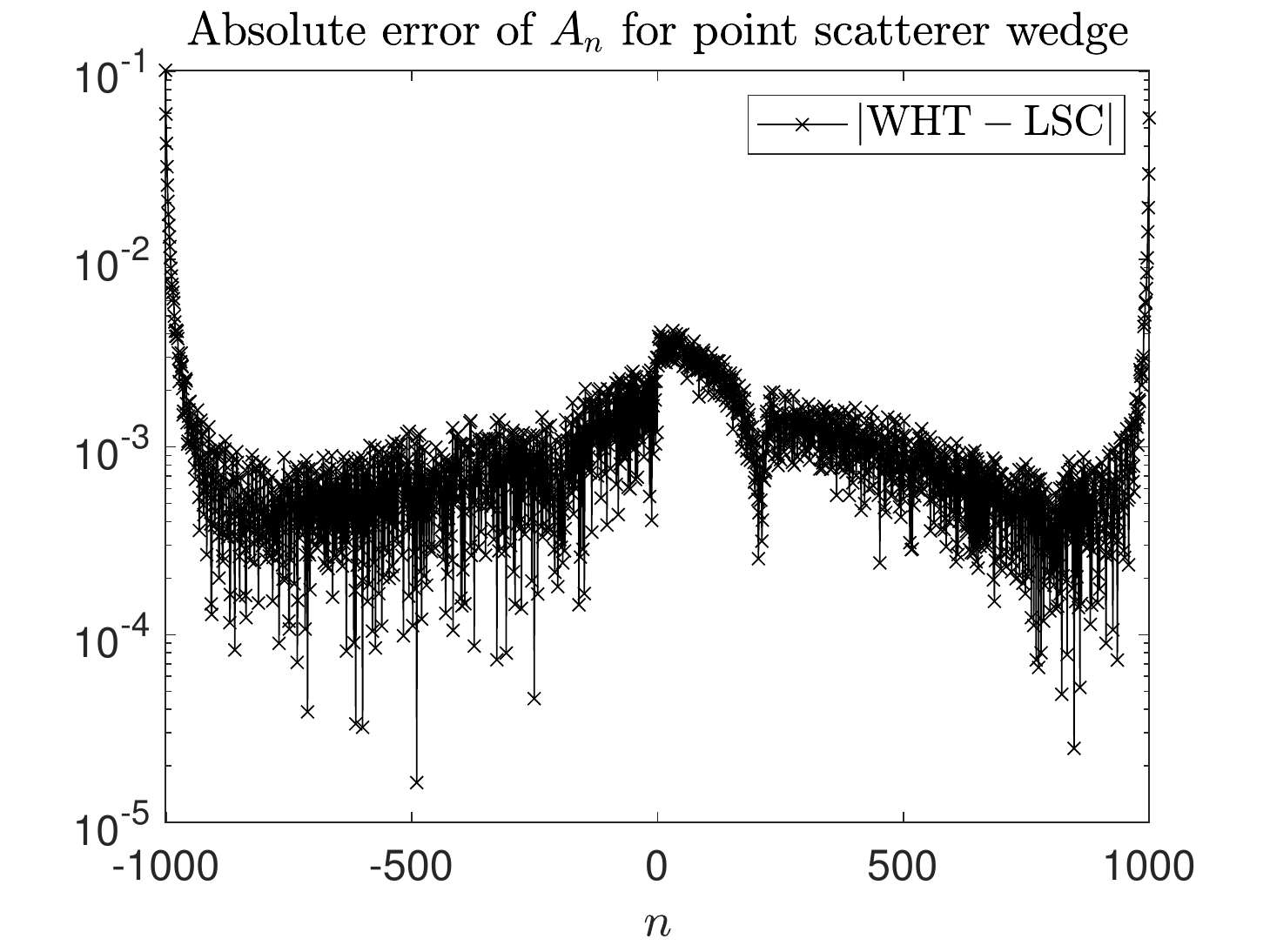}
\caption{Absolute value of the difference between the different methods used to produce the scattering coefficients where the truncation is at $1000$. The index $n$ of the coefficients is ordered such that $A_{-1001},...A_{-1},A_{0},...A_{1000}$ corresponds to $A^{(2)}_{1000},...A^{(2)}_{0},A^{(1)}_{0},...A^{(1)}_{1000}$ respectively. The top row considers an infinite array case with the incident wave having wavenumber $k=5\pi$ and incident angle $\theta_{\textrm{I}}=\frac{\pi}{4}$ (left) or $\theta_{\textrm{I}}=\frac{\pi}{12}$ (right). The bottom row consider the point scatterer wedge case given by Table \ref{table:test-cases} with the incident wave having wavenumber $k=5\pi$ and incident angle $\theta_{\textrm{I}}=0$ (left) or $k=15\pi$ and $\theta_{\textrm{I}}=\frac{\pi}{2}$ (right).}
\label{fig:An-error}
\end{figure}

With this in mind, FIG.\ \ref{fig:An-error} is a collection of plots of the absolute value difference between the two sets of scattering coefficients where the truncation is chosen to be $1000$. In these plots, the index $n$ of the coefficients is ordered such that $A_{-1001},...A_{-1},A_{0},...A_{1000}$ corresponds to $A^{(2)}_{1000},...A^{(2)}_{0},A^{(1)}_{0},...A^{(1)}_{1000}$ respectively. On the top row of FIG.\ \ref{fig:An-error}, we look at the infinite array problem and have an incident wave with wavenumber $k=5\pi$ but two different incident angles; $\theta_{\textrm{I}}=\frac{\pi}{4}$ and $\theta_{\textrm{I}}=\frac{\pi}{12}$ for the left and right side respectively. Note that equivalent plots of the relative error share the same shape and are on similar scale as FIG.\ \ref{fig:An-error}. On the bottom row, we look at the point scatterer wedge given in Table \ref{table:test-cases} and use the same incident waves as in Figure 6 of \cite{HWpaperI} (i.e.\ $k=5\pi$ and $\theta_{\textrm{I}}=0$ for the left side and $k=15\pi$ and $\theta_{\textrm{I}}=\frac{\pi}{2}$ for the right side).

The top row of FIG.\ \ref{fig:An-error} is especially interesting because the comparison with the exact solution allows us to decompose the error between the WH and LSC methods and assess their strengths and weaknesses. One conclusion of this decomposition can be seen in the middle of the plots ($n\approx0$) where the WH seems to be at its weakest and the LSC method at its strongest. This is due to how each method sees the problem as the WH method considers each array separately and adds the interaction between arrays when solved, whereas the LSC method considers each scatterer separately and solves for the individual interactions. The other conclusion of this decomposition can be seen at the truncated ends of the arrays ($n\approx\pm1000$) where the WH is at its strongest and the LSC method at its weakest. This is again due to how the methods see the problem as separate semi-infinite arrays or individual scatterers. It is likely that we will come to similar conclusions for most (if not all) of the potential configurations of this setup. It is also likely that the weak region for the WH method will improve when the semi-infinite arrays are well separated. It is important to note that the overall error does converge to zero as the truncation increases and shape of these error graphs scales with the truncation as well.

\section{Conclusions}
To summarise, we have generalised the WH method used for diffraction by a wedge of point scatterers \cite{HWpaperI}, by considering any number of semi-infinite arrays with an arbitrary set of parameters. The method remains essentially unchanged with the additional benefit of having removed the need for an iterative scheme. The MATLAB scripts we created from this solution are quite versatile and we are able to consider a very wide range of cases, some of which are illustrated in FIG.\ \ref{fig:test-cases}. 

We have also compared the WH method with some numerical approaches such as the LSC method, and FIG.\ \ref{fig:An-error} has shown that the two methods have some good agreement and highlighted the strengths and weaknesses between them. We found that the LSC method is better at modelling the interactions between the scatterers at the ends of the arrays and the WH method is better at modelling the infiniteness of the arrays. Knowing this, one could propose to use a hybrid of the two methods to get accurate coefficients for all scatterers (in other words, use LSC to get the coefficients $A^{(j)}_m$ for small $m$ and WH for large $m$). In theory, this hybrid would have the strengths of both methods but neither of the weaknesses. While finding the optimal $m$ where we should switch methods would be simple for infinite array cases, more general configurations will be more difficult. This is because we will not have an exact solution to decompose the error quantity $|\text{WHT}-\text{LSC}|$ and this optimal $m$ will not be unique.

It is also possible to reformulate the entries of $\bar{\mathcal{M}}^{(j,\ell)}$ \eqref{MSIA-barM} by rewriting the Hankel function in its integral form, evaluating the sum and then approximating the result using the method of steepest descent. This would lead to,
\begin{align}
\bar{\mathcal{M}}^{(j,\ell)}_{nq}\approx\sigma^{(j,\ell)}(n,q)\sqrt{\tfrac{2}{\pi k\Lambda^{(j,\ell)}(n,q)}}e^{ik\Lambda^{(j,\ell)}(n,q)-\frac{i\pi}{4}},
\end{align}
where $\sigma^{(j,\ell)}(n,q)$ is a function of $K^+_j$ which can change depending on the positioning of the scatterers at $\BF{R}^{(j)}_n$ and $\BF{R}^{(\ell)}_q$. The idea here is that an efficient approximation could improve the computation time for large truncations by reducing the computational order with respect to the truncation without sacrificing too much accuracy.

Although we have not explicitly discussed resonance in the current article (see \cite{ArrayResonanceFoldyPaper} for an overview), it is nonetheless of special interest. By exclusively using the methods discussed in this article, we are capable of numerically evaluating cases where inward resonance is occurring, $(ks_j/2\pi)(1+\cos(\alpha_j-\theta_{\textrm{I}}))\in\ZZ$. However, we cannot numerically evaluate outward resonance (and by extension double resonance) cases, $(ks_j/2\pi)(1-\cos(\alpha_j-\theta_{\textrm{I}}))\in\ZZ$, because the associated scattering coefficients will tend to zero but still lead to a non-trivial wave field. Finding a general procedure which can find and extract outward resonant waves is a topic for future work. 

Another interesting avenue to pursue is to find a way to use the solution for the scattering coefficients to identify special features. Examples of this could include; the parameters of Bloch waves in lattices (see the multiple infinite array case given by FIG.\ \ref{fig:test-cases}(f)), trapped modes in waveguides or the constructive/destructive interference caused by waveguides or the acoustic Faraday cage.

\section*{Acknowledgements}
This research was supported by EPSRC grant EP/W018381/1. A.V.K. is supported by a Royal Society Dorothy Hodgkin Research Fellowship and a Dame Kathleen Ollerenshaw Fellowship.
\RED{A.V.K. also acknowledges the support from EU through the H2020-MSCA-RISE2020 project EffectFact, Grant agreement ID: 101008140.} M.A.N. was supported by a David Crighton fellowship at Cambridge University and thanks Nigel Peake for many insightful discussions which greatly refined this article. The authors would like to thank the Isaac Newton Institute for Mathematical Sciences for support and hospitality during the programme \emph{Mathematical theory and applications of multiple wave scattering} when work on this paper was undertaken. This programme was supported by EPSRC grant EP/R014604/1. The authors also give thanks to Stuart Hawkins and David Hewett for providing the source code for the TMAT and LSC methods respectively.

\bibliographystyle{abbrvnat}
\bibliography{HWbib}

\begin{thebibliography}{58}
\providecommand{\natexlab}[1]{#1}
\providecommand{\url}[1]{\texttt{#1}}
\expandafter\ifx\csname urlstyle\endcsname\relax
  \providecommand{\doi}[1]{doi: #1}\else
  \providecommand{\doi}{doi: \begingroup \urlstyle{rm}\Url}\fi

\bibitem[Abrahams and Wickham(1988)]{Abrahams1988}
I.~D. Abrahams and G.~R. Wickham.
\newblock {On the scattering of sound by two semi-infinite parallel staggered
  plates - I. Explicit matrix Wiener-Hopf factorization}.
\newblock \emph{Proc. R. Soc. A}, 420\penalty0 (1858):\penalty0 131--156, 1988.

\bibitem[Abrahams and Wickham(1990{\natexlab{a}})]{Abrahams1988part2}
I.~D. Abrahams and G.~R. Wickham.
\newblock {The scattering of sound by two semi-infinite parallel staggered
  plates. II. Evaluation of the velocity potential for an incident plane wave
  and an incident duct mode}.
\newblock \emph{Proc. R. Soc. A}, 427\penalty0 (1872):\penalty0 139--171,
  1990{\natexlab{a}}.

\bibitem[Abrahams and Wickham(1990{\natexlab{b}})]{Abrahams1990}
I.~D. Abrahams and G.~R. Wickham.
\newblock {Acoustic scattering by two parallel slightly staggered rigid
  plates}.
\newblock \emph{Wave Motion}, 12\penalty0 (3):\penalty0 281--297,
  1990{\natexlab{b}}.

\bibitem[Adams et~al.(2008)Adams, Craster, and
  Guenneau]{AdamsCrasterGuenneau2008}
S.~D. Adams, R.~V. Craster, and S.~Guenneau.
\newblock {Bloch waves in periodic multi-layered acoustic waveguides}.
\newblock \emph{Proc. R. Soc. A}, 464\penalty0 (2098):\penalty0 2669--2692,
  2008.

\bibitem[Baddoo and Ayton(2018)]{Baddoo2018}
P.~J. Baddoo and L.~J. Ayton.
\newblock {Potential flow through a cascade of aerofoils: direct and inverse
  problems}.
\newblock \emph{Proc. R. Soc. A}, 474\penalty0 (2217), 2018.

\bibitem[Baddoo and Ayton(2020)]{Baddoo2020}
P.~J. Baddoo and L.~J. Ayton.
\newblock {An analytic solution for gust-cascade interaction noise including
  effects of realistic aerofoil geometry}.
\newblock \emph{Journal of Fluid Mechanics}, 886, 2020.

\bibitem[Beatson and Greengard(1997)]{BeatsonGreengardFMM}
R.~Beatson and L.~Greengard.
\newblock {A short course on fast multipole methods}.
\newblock In \emph{Wavelets, multilevel methods, and elliptic PDEs}, chapter~1,
  pages 1--37. Oxford University Press, 1997.

\bibitem[{Bonnet-Ben Dhia} and Starling(1994)]{BonnetBenDhia1994}
A.-S. {Bonnet-Ben Dhia} and F.~Starling.
\newblock {Guided waves by electromagnetic gratings and non‐uniqueness
  examples for the diffraction problem}.
\newblock \emph{Math. Methods Appl. Sci.}, 17:\penalty0 305--338, 1994.

\bibitem[{Bonnet-Ben Dhia} et~al.(2016){Bonnet-Ben Dhia}, Fliss, Hazard, and
  Tonnoir]{BonnetBenDhia2016}
A.-S. {Bonnet-Ben Dhia}, S.~Fliss, C.~Hazard, and A.~Tonnoir.
\newblock {A Rellich type theorem for the Helmholtz equation in a conical
  domain}.
\newblock \emph{C. R. Math.}, 354\penalty0 (1):\penalty0 27--32, 2016.

\bibitem[Botten et~al.(2001)Botten, Nicorovici, McPhedran, de~Sterke, and
  Asatryan]{Botten2001}
L.~C. Botten, N.~A. Nicorovici, R.~C. McPhedran, C.~M. de~Sterke, and A.~A.
  Asatryan.
\newblock {Photonic band structure calculations using scattering matrices}.
\newblock \emph{Phys. Rev. E}, 64:\penalty0 046603, 2001.

\bibitem[Chapman et~al.(2015)Chapman, Hewett, and
  Trefethen]{ChapmanHewettTrefethen2015}
S.~J. Chapman, D.~P. Hewett, and L.~N. Trefethen.
\newblock Mathematics of the {F}araday cage.
\newblock \emph{SIAM Rev.}, 57\penalty0 (3):\penalty0 398--417, 2015.

\bibitem[Craster et~al.(2009)Craster, Guenneau, and
  Adams]{CrasterGuenneauAdams2009}
R.~V. Craster, S.~Guenneau, and S.~D. Adams.
\newblock {Mechanism for slow waves near cutoff frequencies in periodic
  waveguides}.
\newblock \emph{Phys. Rev. B}, 79\penalty0 (4):\penalty0 1--5, 2009.

\bibitem[Crutchfield et~al.(2006)Crutchfield, Gimbutas, Greengard, Huang,
  Rokhlin, Yarvin, and Zhao]{Greengardetal2006}
W.~Crutchfield, Z.~Gimbutas, L.~Greengard, J.~Huang, V.~Rokhlin, N.~Yarvin, and
  J.~Zhao.
\newblock {Remarks on the implementation of the wideband FMM for the Helmholtz
  equation in two dimensions}.
\newblock \emph{Contemp. Math.}, 408:\penalty0 99--110, 2006.

\bibitem[Daniele and Zich(2014)]{DanieleZich2014}
V.~G. Daniele and R.~S. Zich.
\newblock \emph{{The Wiener-Hopf Method in Electromagnetics}}.
\newblock Scitech, Edison, New Jersey, 2014.

\bibitem[Foldy(1945)]{Foldy1945}
L.~L. Foldy.
\newblock {The multiple scattering of waves. I. General theory of isotropic
  scattering by randomly distributed scatterers}.
\newblock \emph{Phys. Rev.}, 67\penalty0 (3-4):\penalty0 107--119, 1945.

\bibitem[Ganesh and Hawkins(2017)]{HawkinsGanesh2017}
M.~Ganesh and S.~C. Hawkins.
\newblock {Algorithm 975: TMATROM-A T-matrix reduced order model software}.
\newblock \emph{ACM Trans. Math. Softw.}, 44\penalty0 (1), 2017.

\bibitem[Greengard and Gimbutas(2012)]{FMM_library}
L.~Greengard and Z.~Gimbutas.
\newblock {FMMLIB2D: A matlab toolbox for fast multipole method in two
  dimensions version 1.2}, 2012.
\newblock URL \url{https://cims.nyu.edu/cmcl/cmcl.html}.

\bibitem[Haslinger et~al.(2014)Haslinger, Movchan, Movchan, and
  McPhedran]{Haslinger2014}
S.~G. Haslinger, A.~B. Movchan, N.~V. Movchan, and R.~C. McPhedran.
\newblock {Symmetry and resonant modes in platonic grating stacks}.
\newblock \emph{Waves Random Complex Media}, 24\penalty0 (2):\penalty0
  126--148, 2014.

\bibitem[Haslinger et~al.(2016)Haslinger, Craster, Movchan, Movchan, and
  Jones]{Haslinger2016}
S.~G. Haslinger, R.~V. Craster, A.~B. Movchan, N.~V. Movchan, and I.~S. Jones.
\newblock {Dynamic interfacial trapping of flexural waves in structured
  plates}.
\newblock \emph{Proc. R. Soc. A}, 472\penalty0 (20152186), 2016.

\bibitem[Haslinger et~al.(2018)Haslinger, Jones, Movchan, and
  Movchan]{HaslingerJonesMovchan2018}
S.~G. Haslinger, I.~S. Jones, N.~V. Movchan, and A.~B. Movchan.
\newblock {Localization in semi-infinite herringbone waveguides}.
\newblock \emph{Proc. R. Soc. A}, 474:\penalty0 20170590, 2018.

\bibitem[Hawkins(2023)]{tmatsolver}
S.~C. Hawkins.
\newblock {A T-matrix repository for two- and three-dimensional multiple wave
  scattering simulations in MATLAB}, 2023.
\newblock URL \url{https://github.com/stuart-hawkins/tmatsolver}.

\bibitem[Heins(1948{\natexlab{a}})]{Heins1948part1}
A.~E. Heins.
\newblock {The radiation and transmission properties of a pair of semi-infinite
  parallel plates - I}.
\newblock \emph{Q. Appl. Math.}, 6\penalty0 (2):\penalty0 157--166,
  1948{\natexlab{a}}.

\bibitem[Heins(1948{\natexlab{b}})]{Heins1948part2}
A.~E. Heins.
\newblock {The radiation and transmission properties of a pair of parallel
  plates - II}.
\newblock \emph{Q. Appl. Math.}, 6\penalty0 (3):\penalty0 215--220,
  1948{\natexlab{b}}.

\bibitem[Hewett and Hewitt(2016)]{HewettHewitt2016}
D.~P. Hewett and I.~J. Hewitt.
\newblock Homogenized boundary conditions and resonance effects in {F}araday
  cages.
\newblock \emph{Proc. Roy. Soc. A}, 472\penalty0 (2189), 2016.

\bibitem[Hills and Karp(1965)]{HillsKarp1965}
N.~L. Hills and S.~N. Karp.
\newblock {Semi-Infinite Diffraction Gratings-I}.
\newblock \emph{Comm. Pure Appl. Math}, XVIII:\penalty0 203--233, 1965.

\bibitem[Jones(1986)]{DSJones1986}
D.~S. Jones.
\newblock {Diffraction by three semi-infinite planes}.
\newblock \emph{Proc. R. Soc. A}, 404:\penalty0 299--321, 1986.

\bibitem[Jones et~al.(2017)Jones, Movchan, and Movchan]{JonesMovchan2017}
I.~S. Jones, N.~V. Movchan, and A.~B. Movchan.
\newblock {Blockage and guiding of flexural waves in a semi-infinite double
  grating}.
\newblock \emph{Mathematical Methods in the Applied Sciences}, 40\penalty0
  (9):\penalty0 3265--3282, 2017.

\bibitem[Kirby(2008)]{Kirby2008}
R.~Kirby.
\newblock {Modeling sound propagation in acoustic waveguides using a hybrid
  numerical method}.
\newblock \emph{The Journal of the Acoustical Society of America}, 124\penalty0
  (4):\penalty0 1930--1940, 2008.

\bibitem[Kisil(2018)]{Kisil2018}
A.~V. Kisil.
\newblock {An iterative wiener-hopf method for triangular matrix functions with
  exponential factors}.
\newblock \emph{SIAM J. Appl. Math.}, 78\penalty0 (1):\penalty0 45--62, 2018.

\bibitem[Kisil and Ayton(2018)]{KisilAyton2018}
A.~V. Kisil and L.~J. Ayton.
\newblock {Aerodynamic noise from rigid trailing edges with finite porous
  extensions}.
\newblock \emph{J. Fluid Mech.}, 836:\penalty0 117--144, 2018.

\bibitem[Krynkin and McIver(2009)]{KrynkinMcIver2009}
A.~Krynkin and P.~McIver.
\newblock {Approximations to wave propagation through a lattice of Dirichlet
  scatterers}.
\newblock \emph{Waves Random Complex Media}, 19\penalty0 (2):\penalty0
  347--365, 2009.

\bibitem[Lawrie and Abrahams(2007)]{LawrieAbrahams2007}
J.~B. Lawrie and I.~D. Abrahams.
\newblock {A brief historical perspective of the Wiener-Hopf technique}.
\newblock \emph{J. Eng. Math.}, 59\penalty0 (4):\penalty0 351--358, 2007.

\bibitem[Linton(1998)]{Linton1998}
C.~M. Linton.
\newblock {The {G}reen's function for the two-dimensional {H}elmholtz equation
  in periodic domains}.
\newblock \emph{J. Eng. Math.}, 33\penalty0 (4):\penalty0 377--402, 1998.

\bibitem[Linton(2006)]{Linton2006}
C.~M. Linton.
\newblock {Schl{\"{o}}milch series that arise in diffraction theory and their
  efficient computation}.
\newblock \emph{J Phys. A}, 39\penalty0 (13):\penalty0 3325--3339, 2006.

\bibitem[Linton(2010)]{Linton2010}
C.~M. Linton.
\newblock {Lattice sums for the {H}elmhoitz equation}.
\newblock \emph{SIAM Rev.}, 52\penalty0 (4):\penalty0 630--674, 2010.

\bibitem[Linton and Martin(2004)]{LintonMartin2004}
C.~M. Linton and P.~A. Martin.
\newblock {Semi-Infinite Arrays of Isotropic Point Scatterers. A Unified
  Approach}.
\newblock \emph{SIAM J. Appl. Math.}, 64\penalty0 (3):\penalty0 1035--1056,
  2004.

\bibitem[Lynott et~al.(2019)Lynott, Andrew, Abrahams, Simon, Parnell, and
  Assier]{Lynott_19}
G.~M. Lynott, V.~Andrew, I.~D. Abrahams, M.~J. Simon, W.~J. Parnell, and R.~C.
  Assier.
\newblock Acoustic scattering from a one-dimensional array; tail-end
  asymptotics for efficient evaluation of the quasi-periodic green's function.
\newblock \emph{Wave Motion}, 89:\penalty0 232--244, 2019.

\bibitem[Maierhofer and Peake(2020)]{MaierhoferPeake2020}
G.~Maierhofer and N.~Peake.
\newblock {Wave scattering by an infinite cascade of non-overlapping blades}.
\newblock \emph{J. Sound Vib.}, 481\penalty0 (115418), 2020.

\bibitem[Maierhofer and Peake(2022)]{MaierhoferPeake2022}
G.~Maierhofer and N.~Peake.
\newblock {Acoustic and hydrodynamic power of wave scattering by an infinite
  cascade of plates in mean flow}.
\newblock \emph{J. Sound Vib.}, 520\penalty0 (116564), 2022.

\bibitem[Martin(2006)]{pmartin2006}
P.~A. Martin.
\newblock \emph{{Multiple Scattering: Interaction of Time-Harmonic Waves with N
  obstacles}}.
\newblock Cambridge University Press, Cambridge, 2006.

\bibitem[Martin(2014)]{Martin2014}
P.~A. Martin.
\newblock {On acoustic and electric Faraday cages}.
\newblock \emph{Proc. R. Soc. A}, 470:\penalty0 20140344, 2014.

\bibitem[Maurya and Sharma(2019)]{MauryaSharma2019}
G.~Maurya and B.~L. Sharma.
\newblock {Scattering by two staggered semi-infinite cracks on square lattice:
  an application of asymptotic Wiener–Hopf factorization}.
\newblock \emph{Z. Angew. Math. Phys.}, 70\penalty0 (5):\penalty0 1--21, 2019.

\bibitem[McIver(2007)]{McIver2007}
P.~McIver.
\newblock {Approximations to wave propagation through doubly-periodic arrays of
  scatterers}.
\newblock \emph{Waves Random Complex Media}, 17\penalty0 (4):\penalty0
  439--453, 2007.

\bibitem[Movchan et~al.(2009)Movchan, McPhedran, Movchan, and
  Poulton]{Movchan2009}
N.~V. Movchan, R.~C. McPhedran, A.~B. Movchan, and C.~G. Poulton.
\newblock {Wave scattering by platonic grating stacks}.
\newblock \emph{Proc. R. Soc. A}, 465\penalty0 (2111):\penalty0 3383--3400,
  2009.

\bibitem[Nethercote et~al.(2020{\natexlab{a}})Nethercote, Assier, and
  Abrahams]{LTpaper}
M.~A. Nethercote, R.~C. Assier, and I.~D. Abrahams.
\newblock {High-contrast approximation for penetrable wedge diffraction}.
\newblock \emph{IMA Journal of Applied Mathematics}, 85\penalty0 (3):\penalty0
  421--466, 2020{\natexlab{a}}.

\bibitem[Nethercote et~al.(2020{\natexlab{b}})Nethercote, Assier, and
  Abrahams]{WedgeReview}
M.~A. Nethercote, R.~C. Assier, and I.~D. Abrahams.
\newblock {Analytical methods for perfect wedge diffraction: a review}.
\newblock \emph{Wave Motion}, 93\penalty0 (102479), 2020{\natexlab{b}}.

\bibitem[Nethercote et~al.(2022{\natexlab{a}})Nethercote, Kisil, and
  Assier]{HWpaperI}
M.~A. Nethercote, A.~V. Kisil, and R.~C. Assier.
\newblock {Diffraction of acoustic waves by a wedge of point scatterers}.
\newblock \emph{SIAM J. Appl. Math.}, 82\penalty0 (3):\penalty0 872--898,
  2022{\natexlab{a}}.

\bibitem[Nethercote et~al.(2022{\natexlab{b}})Nethercote, Thompson, Kisil, and
  Assier]{ArrayResonanceFoldyPaper}
M.~A. Nethercote, I.~Thompson, A.~V. Kisil, and R.~C. Assier.
\newblock {Array scattering resonance in the context of Foldy's approximation}.
\newblock \emph{Proceedings of the Royal Society A: Mathematical, Physical and
  Engineering Sciences}, 478\penalty0 (20220604), 2022{\natexlab{b}}.

\bibitem[Peake(1992)]{Peake1992}
N.~Peake.
\newblock {The interaction between a high-frequency gust and a blade row}.
\newblock \emph{J. Fluid Mech.}, 241:\penalty0 261--289, 1992.

\bibitem[Peake and Cooper(2001)]{PeakeCooper2001}
N.~Peake and A.~J. Cooper.
\newblock {Acoustic propagation in ducts with slowly varying elliptic
  cross-section}.
\newblock \emph{J. Sound Vib.}, 243\penalty0 (3):\penalty0 381--401, 2001.

\bibitem[Peake and Kerschen(1997)]{PeakeKerschen1997}
N.~Peake and E.~J. Kerschen.
\newblock {Influence of mean loading on noise generated by the interaction of
  gusts with a flat-plate cascade: upstream radiation}.
\newblock \emph{Journal of Fluid Mechanics}, 347:\penalty0 315--346, 1997.

\bibitem[Peake and Kerschen(2004)]{PeakeKerschen2004}
N.~Peake and E.~J. Kerschen.
\newblock {Influence of mean loading on noise generated by the interaction of
  gusts with a cascade: Downstream radiation}.
\newblock \emph{Journal of Fluid Mechanics}, 515:\penalty0 99--133, 2004.

\bibitem[Rogosin and Mishuris(2016)]{RogosinMishuris2016}
S.~V. Rogosin and G.~S. Mishuris.
\newblock {Constructive methods for factorization of matrix-functions}.
\newblock \emph{IMA J. Appl. Math.}, 81\penalty0 (2):\penalty0 365--391, 2016.

\bibitem[Shanin(1998)]{AVShanin1998}
A.~V. Shanin.
\newblock {Excitation of Waves in a Wedge-Shaped Region}.
\newblock \emph{Acoust. Phys.}, 44\penalty0 (5):\penalty0 592--597, 1998.

\bibitem[Sharma(2015{\natexlab{a}})]{BLSharma_SIAM2015a}
B.~L. Sharma.
\newblock {Diffraction of waves on square lattice by semi-infinite crack}.
\newblock \emph{SIAM J. Appl. Math.}, 75\penalty0 (3):\penalty0 1171--1192,
  2015{\natexlab{a}}.

\bibitem[Sharma(2015{\natexlab{b}})]{BLSharma_WM2015}
B.~L. Sharma.
\newblock {Diffraction of waves on square lattice by semi-infinite rigid
  constraint}.
\newblock \emph{Wave Motion}, 59:\penalty0 52--68, 2015{\natexlab{b}}.

\bibitem[Thompson and Linton(2008)]{ThompsonLinton2008}
I.~Thompson and C.~M. Linton.
\newblock {An interaction theory for scattering by defects in arrays}.
\newblock \emph{SIAM J. Appl. Math.}, 68\penalty0 (6):\penalty0 1783--1806,
  2008.

\bibitem[Thompson et~al.(2008)Thompson, Linton, and
  Porter]{ThompsonLintonPorter2008}
I.~Thompson, C.~M. Linton, and R.~Porter.
\newblock {A new approximation method for scattering by long finite arrays}.
\newblock \emph{Q. J. Mech. Appl. Math}, 61\penalty0 (3):\penalty0 333--352,
  2008.

\end{thebibliography}

\end{document}